\documentclass[10pt]{article}

\usepackage{amssymb,amsthm,amsmath}
\usepackage[usenames,dvipsnames]{xcolor}
\usepackage{mathabx}
\allowdisplaybreaks
\usepackage{color,graphicx,epsfig}
\usepackage{mathrsfs}
\usepackage{enumerate} 
\usepackage{enumitem} 
\usepackage{esvect}

\newcommand{\comment}[1]{}
\newcommand{\cayd}[2]{\overrightarrow{Cay}(#1;#2)}
\newcommand{\cay}[2]{Cay(#1;#2)}
\definecolor{teal}{RGB}{0,128,128}
\definecolor{darkpurple}{RGB}{128,0,128}

\usepackage[titletoc,toc,title]{appendix}

\newtheorem{theorem}{Theorem}[section]

\newtheorem{lemma}[theorem]{Lemma}
\newtheorem{cor}[theorem]{Corollary}

\theoremstyle{definition}
\newtheorem{defn}[theorem]{Definition}
\newtheorem{rem}[theorem]{Remark}
\theoremstyle{definition}

\newtheorem{ex}[theorem]{Example}

\def \cB {{\cal B}}
\def \cC {{\cal C}}

\def \cE {{\cal E}}
\def \cF {{\cal F}}

\def \cL {{\cal L}}

\def \cR {{\cal R}}

\def \cW {{\cal W}}

\def \Z {\mathbb Z}

\def \aa {\mathbf a}
\def \bb {\mathbf b}
\def \cc {\mathbf c}

\def \hh {\mathbf h}

\def \kk {\mathbf k}

\def \u  {\underline}
\newcommand{\z}{zero sum}
\newcommand{\nz}{nonzero sum}
\newcommand{\Nz}{Nonzero sum}

\title{Constructing generalized Heffter arrays via near alternating sign matrices}

\author{
L.\ Mella\footnotemark[1]
,
T.\ Traetta\footnotemark[2]
}
\date{\vspace{-5ex}}
\begin{document}
\maketitle

\footnotetext[1]{
Dip. di Scienze Fisiche, Informatiche, Matematiche, Universit\`{a} degli Studi di Modena e Reggio Emilia, Via Campi 213/A, I-41125 Modena, Italy.
E-mail: lorenzo.mella@unipr.it
}
\footnotetext[2]{DICATAM, Universit\`{a} degli Studi di Brescia, Via Branze 43, 25123 Brescia, Italy. E-mail: tommaso.traetta@unibs.it}

\begin{abstract} 
Let $S$ be a subset of a group $G$ (not necessarily abelian) such that $S\,\cap -S$ is empty or contains only elements of order $2$, 
and let $\mathbf{h}=(h_1,\ldots, h_m)\in \mathbb{N}^m$ 
and $\mathbf{k}=(k_1, \ldots, k_n)\in \mathbb{N}^n$.
A \emph{generalized Heffter array} GHA$^{\lambda}_S(m, n; \mathbf{h}, \mathbf{k})$ over $G$ is an $m\times n$ matrix $A=(a_{ij})$ such that: the $i$-th row (resp. $j$-th column) of $A$ contains exactly $h_i$ (resp. $k_j$) nonzero elements, 
and the list $\{a_{ij}, -a_{ij}\mid a_{ij}\neq 0\}$ equals $\lambda$ times the set $S\,\cup\, -S$. 
We speak of a \z{} (resp. \nz{}) GHA if 
each row and each column of $A$ sums to zero (resp. a nonzero element),
with respect to some ordering. 

In this paper, we use \emph{near alternating sign matrices} to build 
both zero and \nz{} GHAs, over cyclic groups,
having the further strong property of being simple. 
In particular, we construct zero sum and simple GHAs whose row and column weights are congruent to $0$ modulo $4$. 
This result also provides the first infinite family of simple (classic) Heffter arrays to be rectangular ($m\neq n$) and with less than $n$ nonzero entries in each row.
Furthermore, 
we build \nz{} GHA$^{\lambda}_S(m, n; \mathbf{h}, \mathbf{k})$  over an arbitrary group $G$ whenever $S$ contains enough noninvolutions, thus extending previous nonconstructive results where  $\pm S = G\setminus H$ for some subgroup $H$~of~$G$.

Finally, we describe how GHAs can be used to build orthogonal decompositions and biembeddings of Cayley graphs (over groups not necessarily abelian) onto orientable surfaces.
\end{abstract}

\section{Introduction}
In \cite{A}, Archdeacon introduced Heffter arrays as a tool to construct current graphs, orthogonal cycle systems and biembeddings of complete graphs on surfaces. Since then, such arrays have been subject to extensive study and various generalizations, see for instance \cite{ABD, ADDY, BU, BP, BCDY,CDY, CDF, CDFP, CMPP20, CPEJC, DW}. In the following, we start by introducing a generalized class of Heffter-type arrays, thus unifying different terminologies used in quite a few papers concerning Heffter arrays and some variants, for which we refer the reader to the survey \cite{PD23}.

Let $G$ be an additive group, not necessarily abelian, and
choose a map $\|\cdot\|:G \rightarrow G, a \mapsto \|a\|$
such that $\|a\|=\| {-a} \|\in\{\pm a\}$, for every $a\in G$.
We refer to $\|a\|$ as the \emph{absolute value} of $a\in G$.
Given a multiset (or list) $H$ of $G$, we define the multiset
$\|H\| = \{\|h\| \mid h\in H\}$, while $H^+$ represents the underlying set of $\|H\|$ minus the zero element. In particular, set $\mathbb{Z}^+=\mathbb{N}=\{1,2,\ldots\}$. For a set $S\subseteq G$, we denote by $I(S)$ the set of all involutions (i.e., elements of order $2$) in $S$.

For a multiset $M$ and a positive integer $\lambda$, we denote by 
${}^\lambda M$ the multiset union of $\lambda$ copies of $M$.

Given an $m\times n$ matrix $A$ with entries from $G$, 
the \emph{row-weight} of $A$ is the sequence $w_r(A)=(h_1, \ldots, h_m)$ where $h_i$ is the number of nonzero entries in the $i$-th row of $A$; the \emph{column-weight} $w_c(A)$ of $A$ 
is the row-weight of the transpose of $A$. 
We denote by $\mathcal{E}(A)$ the multiset of nonzero entries of $A$. 

\begin{defn}\label{GHA}
Let $S\subseteq G^+$, $(m,n,\lambda)\in \mathbb{N}^3$,  and let $\mathbf{h}=(h_1, \ldots, h_m)$ and $\mathbf{k}=(k_1, \ldots, k_n)$,
with $1\leq h_i\leq n$ and $1\leq k_j\leq m$ for each $i$ and $j$.
A \emph{generalized Heffter array} GHA$^{\lambda}_S(m, n; \mathbf{h}, \mathbf{k})$ (briefly, GHA) over $G$
is an $m\times n$ matrix $H$, with entries from $G$, satisfying the following properties:
\begin{enumerate}
   \item $\;^{2}\|\cE(H)\| = 
   \;^{2\lambda} (S\setminus I(S)) \,\cup\, ^{\lambda}I(S)$, and  
  \item the $i$-th row (resp. $j$-th column) of $A$ contains exactly $h_i$ (resp. $k_j$) nonzero elements, that is, $(w_r(H), w_c(H)) = (\hh, \kk)$.
\end{enumerate}
We drop the parameter $\lambda$ (representing the \emph{multiplicity} of a GHA) when it is equal to 1.
We speak of a \emph{uniform} GHA whenever the weight sequences $\hh$ and $\kk$ are constant, that is, $\hh = (h,\ldots, h)$ and $\kk = (k, \ldots, k)$; in this case, we write GHA$^{\lambda}_S(m, n; h, k)$. 
\end{defn}
It follows that for a GHA$^{\lambda}_S(m, n; \mathbf{h}, \mathbf{k})$ to exist, we necessarily have that
\begin{equation}\label{GHA:nec}
\begin{aligned} 
 \text{$\lambda$ must be even when $I(S)\neq \varnothing$, and}\\
 \lambda|S| - \textstyle{\frac{\lambda}{2}} |I(S)|=h_1 + \ldots + h_m = k_1 + \ldots + k _n.
\end{aligned}
\end{equation}
Clearly, changing the signs of some entries of a GHA produces another GHA with the same parameters.

We recall the following result by Gale and Ryser \cite[Theorem 7.7.4]{DJ} on the existence of matrices over $\Z_2$ with given row and column weights.
\begin{theorem}\label{GR}
	  Let $\hh=(h_1, \ldots, h_m)$ and $\kk=(k_1, \ldots,k_n)$ be two sequences of positive integers such that $\textstyle{\sum_{i=1}^{m} h_i=\sum_{j=1}^{n} k_j}$. 
	  There exists an $m\times n$ matrix 
	  $A$ over $\mathbb{Z}_2$ with $w_r(A)=\hh$ and $w_c(A)=\kk$ if~and~only~if 
  \begin{equation}\label{GHA:nec2}
    \text{$\textstyle \sum_{i=1}^m \min (h_i,u) \geq \sum_{j=1}^u k'_j$, for every $1\leq u\leq n$}.
  \end{equation}
  where $(k'_1, \ldots, k'_n)$ is the decreasing reordering of $\kk$.
\end{theorem}

Replacing each nonzero entry of a 
GHA$^{\lambda}_S(m, n; \mathbf{h}, \mathbf{k})$
with $1$ produces a matrix over $\Z_2=\{0,1\}$ whose row and column weights are 
$\mathbf{h}$ and $\mathbf{k}$. Therefore, \eqref{GHA:nec2} provides another necessary condition for the existence of a GHA. We will refer to \eqref{GHA:nec} and \eqref{GHA:nec2} as \emph{the necessary conditions} for the existence of a GHA. 

As described in Sections \ref{sec:decs} and \ref{sec:biembeds}, GHAs can be used to construct orthogonal path or cycle decompositions and biembeddings of Cayley graphs onto orientable surfaces. The structural properties of these decompositions and biembeddings depend on the sum of the entries in each row and column of a GHA, with respect to a given ordering. Indeed, it is needless to recall that the sum of elements in a non abelian group heavily depends on the chosen ordering. Given a sequence 
$\aa=(a_1, \ldots, a_n)\in G^n$, we denote by $s_0(\aa)=\sum_{i=1}^n a_i$ the sum of 
$\aa$.

Let $A$ be an $m\times n$ matrix over a group $G$, and denote by $A_i$ (resp. $A^j$) its $i$-th row (resp. $j$-th column). An \emph{ordering} -- of the nonzero entries -- of $A_i$ (resp. $A^j$) is a sequence $\omega_{r,i}$ (resp. $\omega_{c,j}$) whose underlying set is $\mathcal{E}(A_i)$ (resp. $\mathcal{E}(A^j)$). 
We say that $\omega_{r,i}$ (resp. $\omega_{c,j}$) is a \emph{natural} ordering 
if it is obtained by cyclically sequencing from left to right (resp. from top to bottom) the nonzero entries of $A_i$ (resp. $A^j$) starting from a given element.

We refer to $\omega_{r} = (\omega_{r,1},\ldots, \omega_{r,m})$ and
$\omega_{c} = (\omega_{c,1},\ldots, \omega_{c,n})$ as orderings of the rows and of the columns of $A$, respectively; we also call 
$\omega=(\omega_{r}, \omega_{c})$ an ordering of $A$. We refer to $(A,\omega)$ as an \emph{ordered} matrix (with respect to the ordering $\omega$).
We say that $\omega_{r}$ (resp. $\omega_{c}$) is a natural ordering of the rows (resp. columns) of $A$ if each $\omega_{r,i}$ (resp. $\omega_{c,j}$) is natural. Similarly,
$\omega=(\omega_{r}, \omega_{c})$ is natural, if both $\omega_{r}$ and  $\omega_{c}$ are so.

\begin{defn}\label{GHA:zeronon}
  An ordered $m\times n$ matrix $(A,\omega)$ over a group $G$ is 
  \begin{enumerate}
    \item \z{}, if $s_0(\omega_{r,i})=0=s_0(\omega_{c,j})$ 
    for every 
    $1\leq i\leq m$ and $1\leq j\leq n$;
    \item \nz{}, if $s_0(\omega_{r,i})\neq0\neq s_0(\omega_{c,j})$ 
    for every 
    $1\leq i\leq m$ and $1\leq j\leq n$.
  \end{enumerate}
  The matrix $A$ is \z{} (resp. \nz{}) if it is so with respect to some ordering.
\end{defn}
 Note that when $G$ is abelian, the matrix $A$ cannot be at the same time \z{} and \nz; in this case, the property of being
 a \z{} or \nz{} matrix is independent of the chosen ordering. 
 Furthermore, given an arbitrary group $G$ and considering that a cyclic shift of a sequence $\aa$ of $G$ turns it into another sequence whose sum is a conjugate\footnote{We recall that the conjugate of an element $g\in G$ is any group element of the form $-x+g+x$, for some $x\in G$.} of $s_0(\aa)$, it follows that a \z{} matrix over $G$ with respect to a natural ordering continues to be \z{} with respect to every natural ordering.
 
 We point out that the classic concept of a Heffter array, introduced in \cite{A}, is equivalent to a \z{} uniform GHA$^{\lambda}_S(m, n; h, k)$ over the cyclic group $G=\Z_{v}$, where $\lambda=1$, $v = 2hm+1=2kn+1$ and $S=\{1, \ldots, hm\}$.
We notice also that for different values of the parameters of a GHA, we can obtain some variants to the original Heffter arrays for which we refer the reader to the recent survey
in \cite{PD23}.

\begin{defn} A naturally ordered generalized Heffter array (NGHA) is a GHA endowed with a natural ordering.
\end{defn}

In Section \ref{algo}, we construct \nz{} NGHA$^{\lambda}_S(m, n; \mathbf{h}, \mathbf{k})$ over an arbitrary group whenever the necessary conditions hold and $S$ contains enough noninvolutions (elements of order greater than 2), thus extending a result in \cite{CDF}, which is however nonconstructive.
More precisely, we prove the following.

\begin{theorem}\label{thm:main_gh}
Let $G$ be an arbitrary group and let $S\subseteq G^+$ such that
\begin{equation}
 |S\setminus I(S)| \geq m+n-1.
\end{equation}
  There exists a \nz{} NGHA$^\lambda_S(m, n; \mathbf{h}, \mathbf{k})$ over $G$ if and only if the necessary conditions \eqref{GHA:nec} and \eqref{GHA:nec2} hold.
\end{theorem}
 
Further structural properties of decompositions and biembeddings originating from a generalized Heffter array depend on the partial sums of its rows and columns.
To this matter, it is fundamental the concept of a simple sequence over a group $G$: $\aa=(a_1, \ldots, a_n)\in G^n$ is said to be \emph{simple} if
all its partial sums $a_1, a_1+a_2, a_1+a_2+a_3, \ldots, s_0(\aa)$ are pairwise distinct and nonzero, except possibly for $s_0(\aa)=0$.
An ordering of an $m\times n$ matrix $A$ over a group $G$ is said to be \textit{simple}
if the corresponding orderings of the rows,
$\omega_{r,1}, \ldots, \omega_{r,m}$, 
and of the columns, $\omega_{c,1}, \ldots, \omega_{c,n}$, of $A$ are all simple.
  
\begin{defn}\label{def:simpleGHA} If $\omega$ is a simple ordering of an $m\times n$ matrix $A$ over a group $G$,
we say that $(A,\omega)$ is simple. Furthermore,
  \begin{enumerate}
    \item $A$ is simple if it has a simple ordering;
    \item $A$ is \emph{naturally simple} if it has a simple natural ordering, that is, $(A,\omega)$ is simple
    and $\omega$ is a natural ordering of $A$.
  \end{enumerate}
\end{defn}
We point out that classic (nonzero) Heffter arrays with a simple natural ordering were first studied in
\cite{CMPPHeffter} under the name of \emph{globally simple} Heffter arrays.
We notice that unlike for the property of being zero sum, the  simplicity of a sequence is not, in general, preserved under the action of a cyclic permutation of its elements.

In Section \ref{globallysimple:nonzero sum}, we build \nz{} GHAs over the cyclic group that are, in addition, naturally simple, thus extending some previous constructions contained in \cite{CDFP, MP23} concerning the particular case where $S=(\Z_{uw}\setminus u\Z_{uw})^+$,  and both $\hh$ and $\kk$ are constant. Among other things, we prove the following.

\begin{theorem}\label{thm:gl_simple_GH}
There exists a \nz{} and simple NGHA$_S(m, n; \hh, \kk)$ over $\Z_v$ 
 whenever the necessary condition \eqref{GHA:nec} holds 
 $($that is, $S\subseteq [1,\lfloor \frac{v-1}{2}\rfloor]$ 
 and $|S|= s_0(\hh) = s_0(\kk))$, 
 and either
 \begin{enumerate}
 \item $\hh$ and $\kk$ are constant, or 
 \item $m$ and $n$ are even, 
 $\hh = 2\cdot(h_1, h_1, \ldots, h_m, h_m)$, 
 $\kk = 2\cdot(k_1, k_1, \ldots,$ $k_n, k_n)$, and the necessary condition \eqref{GHA:nec2} holds.
 \end{enumerate}      
 \end{theorem}
 
 Our methods also applies in the much harder case of constructing zero sum and simple GHAs. In Section \ref{globallysimple:zero sum}, among other things, we obtain the following result.
\begin{theorem}\label{fromNASMtoGHA2:cor2} 
 Let $v=(2d+1)u \equiv u \;(\mathrm{mod}\; 16)$, where $u=1$ or $u\equiv 0 \;(\mathrm{mod}\; 4)$, and let $U$ be the subgroup of $\Z_v$ of order $u$. If the necessary conditions \eqref{GHA:nec} and \eqref{GHA:nec2} hold, then there exist
 \begin{enumerate}
   \item 
 a zero sum and simple GHA$_{(\Z_v\setminus U)^+}(4m,$ $2n; (4\hh, 4\hh), 4\kk)$ with
  \[\hh=(h_1, h_1, \ldots, h_m, h_m)\;\;\;\text{and}\;\;\;
    \kk=(k_1, k_1, \ldots, k_n, k_n);
  \] 
  \item a zero sum and simple 
GHA$_{(\Z_v\setminus U)^+}(2m,n; 4h, 4k)$.
  \end{enumerate}
\end{theorem}
We point out that Theorem \ref{fromNASMtoGHA2:cor2}.(2) constructs simple (classic) Heffter arrays in cases whose existence was previously unknown, that is, when the arrays are rectangular and with row and column weights less than $m$ and $n$, respectively.
 
The above results rely on the concept of a \emph{near alternating sign matrix}, introduced in Section \ref{NASM}, where they are built in the uniform case (Theorem \ref{thm:ASM}) and in some nonuniform case (Theorem \ref{thm:nonuniform_2}).

In Section \ref{sec:decs} we point out that a GHA (and more generally a matrix) over an arbitrary group $G$ can be used to build \emph{$G$-regular orthogonal decompositions} of Cayley graphs, directed or not, into walks. These walks may be cycles or paths and this depends on the property of some row or column of $A$ to be simple and either \z{} or not. 
In particular, Theorem \ref{thm:ortodecs_4} constucts pairs of orthogonal cycle decompositions, of complete equipartite graphs, whose cycle lengths are congruent to $0$ modulo $4$.
Furthermore, Theorem \ref{thm:ortodecs_3} constructs $\Z_v$-regular orthogonal path decompositions of any Cayley graph over $\Z_v$ (also called a circulant graph) for almost all possible lengths of the paths. 
Finally, in Section \ref{sec:biembeds}, we show that (as for classic Heffter arrays) GHAs can be used to
build biembeddings of Cayley graphs onto orientable surfaces 
(Theorem \ref{prop:biembed_cay}). 
We point out that the embedded Cayley graph does not need be defined over an abelian group, as always supposed in the context of classic Heffter arrays.

The following section collects the terminology and notation used throughout the article, except for some concepts
already defined in the introduction.

\section{Terminology and notation}
Let $G$ be an additive group, not necessarily abelian, and let $\|\cdot\|:G \rightarrow G, a \mapsto \|a\|$
 be a map such that $\|a\|=\| {-a} \|\in\{\pm a\}$, for every $a\in G$.
Recall that $\|a\|$ is referred to as the \emph{absolute value} of $a\in G$.
This concept extends naturally to multisets of $G$, and coordinate-wise to $n$-tuples or arrays with entries from $G$. In particular, we recall that given a multiset $H$ of $G$, 
$H^+$ represents the underlying set of $\|H\|= \{\|h\| \mid h\in H\}$ minus the zero element. 

Letting $v\in \mathbb{N}\,\cup\, \{\infty\}$, we denote by $\Z_v$ the ring of integers modulo $v$ when $v \neq \infty$, otherwise $\Z_\infty=\Z$.
Let $\Z_v^+ = \mathbb{N}$ when $v=\infty$, 
otherwise $\Z_v^+ = \{1,\ldots, \lfloor v/2 \rfloor\}$. 
We extend to $\Z\cup\{\infty\}$ the natural total ordering of $\Z$
by considering $\infty$ as its maximal element.
In $\Z\cup\{\infty\}$, we denote by 
$[u,v]=\{u, u+1, \ldots, v\}$ the interval of integers from $u$ to $v$ whenever $u\leq v$, otherwise $[u,v]=\varnothing$.

Given a set $S\subseteq G$, we denote by $S^{m,n}$ the set of 
$m\times n$ matrices with entries from $S$, and let $S^n = S^{1,n}$. 
Letting $A=(a_{ij})\in G^{m,n}$, we denote by
\begin{enumerate}
  \item $\cE(A)=\{a_{ij}\mid a_{ij}\neq 0\}$ the list (multiset) of nonzero entries of $A$;
  \item $skel(A)=\{(i,j)\mid a_{ij}\neq 0, 1\leq i\leq m, 1\leq j\leq n\}$ the 
  \emph{skeleton of $A$}, that is, the set of positions $(i,j)$ whose relative entry in $A$ is nonzero.
\end{enumerate}
We will always consider $skel(A)$ lexicographically ordered.

\subsection{Sequences over an arbitrary group $G$}
Given an integer $q$, an element $a$ of a group $G$, and two finite sequences $\mathbf{a}=(a_1, \ldots, a_n)$ and
$\mathbf{b}=(b_1, \ldots, b_t)$ over $G$
and length $n$ and $t$, respectively, we use the following terminology and notation.
\begin{enumerate}
  \item $\u{a}=(a, \ldots, a)$. The length of 
  $\u{a}$ will be always clear from the context.
  \item $(\mathbf{a}, \mathbf{b})=(a_1, \ldots, a_n, b_1, \ldots, b_t)$, $q\mathbf{a} = (qa_1, \ldots, qa_n)$, and $-\mathbf{a} = (-1)\mathbf{a}$. 
  \item $s_{i,j}(\mathbf{a}) = \sum_{\ell=i}^j a_\ell$ is a \emph{run} of $\aa$, for every $1\leq i \leq j \leq n$, and we call it  
  \textit{proper} when $(i,j)\neq (1,n)$.
  \item $s_j(\mathbf{a}) = s_{1,j}(\mathbf{a})$ is the \emph{$j$-th partial sum} of  $\mathbf{a}$, for $1\leq j\leq n$.
        For our convenience, we set $s_0(\mathbf{a})=s_n(\mathbf{a})$.
  \item 
$\mathbf{a}$ is called \emph{simple} if all proper runs of $\mathbf{a}$ are different from zero; this is equivalent to saying that $0\neq s_i(\mathbf{a}) \neq s_j(\mathbf{a})$ for all $1 \leq i < j \leq n$. 
  \item $\mathbf{a}$ is called a \emph{\z{} sequence} (resp. \emph{\nz{} sequence}) 
        if $s_0(\aa)=0$ (resp. $s_0(\aa) \neq 0$).
  \item The \emph{alternated forms} of $\mathbf{a}$ 
are the sequences $\mathbf{a}^{\pm} = (-a_1, \ldots,$ $(-1)^ia_i,$ $\ldots, (-1)^na_n)$ and $-\mathbf{a}^{\pm}$.
\end{enumerate}

\begin{rem}\label{plus_minus_a:rem} 
If $G$ is abelian, given an integer $q$ and a sequence $\mathbf{a}=(a_1, \ldots,$ $a_n)\in G^n$, 
we have that $s_i(q\mathbf{a}) = qs_i(\mathbf{a})$, hence  $s_i(-\mathbf{a}) = -s_i(\mathbf{a})$, for  $1\leq i\leq n$. 
\end{rem}

\subsection{Matrices over an arbitrary group $G$}
Given an $m\times n$ matrix $A=(a_{ij})$ with entries from an arbitrary group $G$,  not necessarily abelian, we use the following terminology and notation. We denote by $A_i$ and $A^j$ the $i$-th row and the $j$-th column of $A$, respectively. The \emph{reduced $i$-th row} of $A$, denoted by $A_{(i)}$, is the left-to-right sequence of nonzero elements in
$A_i$. Similarly, the \emph{reduced $j$-th column} of $A$, denoted by $A^{(j)}$, is 
  the top-to-bottom sequence of nonzero elements in $A^j$.
  
  The \emph{row-weight} of $A$ is the sequence $w_r(A) = (w_1(A), \ldots, w_m(A))$
  where each $w_i(A)$ is the length of 
 $A_{(i)}$, that is, the number of nonzero entries in the $i$-th row of $A$.
  The \emph{column-weight} of $A$ is the sequence $w_c(A) = (w^1(A), \ldots, w^n(A))$ 
  where each $w^j(A)$ is the length of 
  $A^{(j)}$, that is, the number of nonzero entries in the $j$-th column of $A$; clearly, $w_c(A) = w_r(A^t)$, where $A^t$ denotes the transpose of $A$. 
  The \emph{weight} of $A$ is the number $w(A)$
  of nonzero entries of $A$. Clearly,
  \[
  w(A)= \sum_{i=1}^{m} w_i(A) = \sum_{j=1}^n w^j(A)= |\cE(A)|.
  \]
  
  The \emph{nonzero position matrix} associated to $A$ is the $m\times n$ array $A^*=(a^*_{ij})$ where $a^*_{ij}=p$ if $(i,j)$ is the $p$-th element in $skel(A)$ endowed with the lexicographic order, otherwise $a^*_{ij}=0$.  
In other words, $a^*_{ij}$ counts the number of nonzero entries of $A$ from $a_{11}$ up to $a_{ij}\neq 0$, using the lexicographic order over the indices $(i,j)$.

  Given a map $f: G\rightarrow G$, 
  let $f(A)=(f(a_{ij}))$ denote the $m\times n$ matrix obtained
  by applying $f$ element-wise on $A$.
  
  Letting $B=(b_{ij})\in \Z^{m,n}$, the Hadamard product of
  $B$ and $A$ is the $m\times n$ matrix $B\circ A=(b_{ij}a_{ij})$. 

\section{Near alternating sign matrices}
\label{NASM}

In this section, we 
construct near alternating sign matrices, which we use in the following to build simple NGHAs.
\begin{defn}
An $m\times n$ matrix $A$ with entries from $\{0,\pm1\}$ is called a \emph{near alternating sign matrix} (NASM) if the nonzero entries of each row and each column alternate.

Letting $\mathbf{h}=w_r(A)$ and $\mathbf{k}=w_c(A)$, we say that
$A$ is a NASM$(m,n; \mathbf{h}, \mathbf{k})$
or simply a NASM$(m,n)$. Finally, if $\mathbf{h}=(h, \ldots, h)$ and $\mathbf{k}=(k, \ldots, k)$, we say that $A$ is uniform and write NASM$(m,n; h,k)$. 
\end{defn}

Near alternating sign matrices were first considered in
\cite{BruKim}, although the terminology was suggested in a private communication by R. Brualdi. Here, we are interested in building them with given row-weight and column-weight sequences.

\begin{rem}\label{NASM:nec}
  Since the absolute value of a NASM$(m,n;\hh,\kk)$ is a $\{0,1\}$-matrix, 
  by Theorem \ref{GR} it follows that $\hh$ and $\kk$ must satisfy condition \eqref{GHA:nec2}.
\end{rem}

We now define the concept of frame of a NASM, a parameter that will enable us to join two or more suitable NASMs to obtain larger NASMs.

\begin{defn}
  Let $A$ be a NASM$(m,n)$. The \emph{frame of $A$} is the quadruple 
  $\varphi(A)=
  \left(A^\leftarrow, A^\rightarrow, A^\uparrow, A^\downarrow\right)$    
  defined as follows:
  \begin{enumerate}
    \item $A^\leftarrow = (a_1, \ldots, a_m)$ where each $a_i$ is the first nonzero entry of the $i$-th row of $A$;
    \item $A^\rightarrow = (a_1, \ldots, a_m)$ where each $a_i$ is the last  nonzero entry of the $i$-th row of $A$;
    \item $A^\uparrow = (A^t)^\leftarrow$ and $A^\downarrow = (A^t)^\rightarrow$, where $A^t$ is the transpose of $A$.
  \end{enumerate} 
\end{defn}

The following two lemmas are straightforward.

\begin{lemma} \label{AAAA}
An $f\times g$ block matrix $B=(A_{ij})$ is a NASM whenever the following types of submatrices of $B$ are all NASM
\[
A_{ij},\qquad[A_{ij}\quad A_{i,j+1}],\qquad
\left[
\begin{array}{c}
A_{ij} \\
A_{i+1,j}
\end{array}
\right].
\]
\end{lemma}

\begin{lemma} \label{AA}
Let $A_i$ be a NASM$(m_i,n_i; \mathbf{h}_i, \mathbf{k}_i)$, for $i=1,2$.
\begin{enumerate} 
\item $[A_1\;\;A_2]$ is a NASM$(m_1,n_1+n_2; \mathbf{h}_1 + \mathbf{h}_2, (\mathbf{k}_1, \mathbf{k}_2))$ if and only if 
 $m_1=m_2$ and $A_1^\rightarrow = -A_2^\leftarrow$. 
\item $
\left[
\begin{array}{c}
A_1 \\
A_2
\end{array}
\right]
$
is a NASM$(m_1+m_2,n_1; (\mathbf{h}_1, \mathbf{h}_2), \mathbf{k}_1+\mathbf{k}_2)$
if and only if $n_1=n_2$ and $A_1^\downarrow = -A_2^\uparrow$.
\end{enumerate}
\end{lemma}

We now show how to build uniform NASMs.

\begin{theorem}\label{thm:ASM}
  There exists a NASM$(m, n; h, k)$ if and only if $mh=nk$.
\end{theorem}
\begin{proof} 

If a NASM$(m, n; h, k)$ exists, then necessarily $mh=nk$. To prove sufficiency, 
let $f=\gcd(m,k)$ and $g=\gcd(n,h)$. Since $mh=nk$, we have that
$\frac{m}{f} \frac{h}{g}= \frac{n}{g}\frac{k}{f}$. Since
$\gcd(\frac{m}{f}, \frac{k}{f}) = 1 = 
 \gcd(\frac{h}{g}, \frac{n}{g})$, we can write
 $\frac{m}{f} = \frac{n}{g} =\ell$, hence $\frac{h}{g} = \frac{k}{f} =  d$.

We start by constructing a NASM$(\ell,\ell; d,d)$.
Let $A=(a_{ij})$ be the $\ell\times \ell$ matrix over $\{0,\pm1\}$ defined as follows:
\[
a_{ij}=
\begin{cases}
  (-1)^{i+j} & \text{if $j\leq i \leq \min(\ell, j+d-1)$, or} \\
             & \text{\;\; $1\leq i\leq j+d-\ell-1$ and $\ell-d$ is even},\\
  (-1)^{i+j+1} & \text{if $1\leq i\leq j+d-\ell-1$ and $\ell-d$ is odd},\\
  0   & \text{ otherwise}.
\end{cases}
\]
One can easily check that each reduced row (resp. reduced column) of $A$ is alternating
of length $d$, that is, $A$ is a NASM$(\ell,\ell;d,d)$. Furthermore,
\begin{equation}\label{A}
  A^\rightarrow = (-1)^{d+1}A^\leftarrow
    \qquad\text{and}\qquad
  A^\downarrow = (-1)^{d+1} A^\uparrow.
\end{equation}

Now we are going to show that the matrix obtained by suitable repetitions of $A$ or $-A$ is a NASM$(m,n;h,k)$. Let $B=(b_{ij})$ be the $f\times g$ block matrix defined as follows:
\[
b_{ij}=
\begin{cases}
  A & \text{if $d$ is even},\\
  (-1)^{i+j} A & \text{if $d$ is odd},\\
\end{cases}
\]
for $1\leq i\leq f$ and $1\leq j\leq g$. Clearly, $B$ is a matrix with $m=f\ell$ rows and 
$n=g\ell$ columns. Also, its 
reduced rows (resp. columns) have each length $h=gd$ (resp. $k=fd$).
Finally, by \eqref{A} and Lemma \ref{AA}, we have that
$\pm[A\quad (-1)^{d}A]$ and 
$
\pm\left[
\begin{array}{c}
A \\
(-1)^{d}A
\end{array}
\right]
$
are NASMs. It then follows by Lemma \ref{AAAA} that
$B$ is a NASM$(m,n; h, k)$.
\end{proof}

\begin{ex}\label{ex:NASM}
In this example we construct a NASM$(6,9;6,4)$ by using the procedure described in the proof of Theorem \ref{thm:ASM}. We have that:
\[
\begin{aligned}
f &= \gcd(m,k) = 2, 		 \qquad g = \gcd(n,h) = 3, \\
\ell &= \frac{m}{f} = 3, \qquad d = \frac{h}{g} = 2. \\  
\end{aligned}
\] 
We can then construct the following NASM$(3, 3; 2, 2)$:
\[
A = \begin{array}{|r|r|r|}\hline
1 & 0 & -1 \\ \hline
-1 & 1 & 0 \\ \hline
0 & -1 & 1\\ \hline
\end{array}
\]
Since $d$ is even, we build the following array $B$, that is a NASM$(6,9;6,4)$:
\[
\begin{aligned}
B&= \begin{array}{|r|r|r|}\hline
A & A & A  \\ \hline
A & A & A  \\ \hline
\end{array} 
& =  \begin{footnotesize}
\begin{array}{||r|r|r||r|r|r||r|r|r||}\hline \hline
1 & 0 & -1 &1 & 0 & -1 &1 & 0 & -1   \\ \hline 
-1 & 1 & 0 & -1 & 1 & 0 & -1 & 1 & 0  \\ \hline
0 & -1 & 1 & 0 & -1 & 1 & 0 & -1 & 1  \\ \hline \hline
1 & 0 & -1 &1 & 0 & -1 &1 & 0 & -1  \\ \hline 
-1 & 1 & 0 & -1 & 1 & 0 & -1 & 1 & 0\\ \hline
0 & -1 & 1 & 0 & -1 & 1 & 0 & -1 & 1  \\ \hline \hline
\end{array}
\end{footnotesize}
\end{aligned}
\]
In order to show the complete procedure of the proof of Theorem \ref{thm:ASM}, we also construct a NASM$(12,16;4,3)$. As before, we compute the parameters:
\[
\begin{aligned}
f &= \gcd(m,k) = 3, 		 \qquad g = \gcd(n,h) = 4, \\
\ell &= \frac{m}{f} = 4, \qquad d = \frac{h}{g} = 1. \\  
\end{aligned}
\]
We take as a NASM$(4,4;1,1)$ the $4\times 4$ identity matrix, denoted by $I$. Since $d$ is odd, we construct the following block matrix, that is a NASM$(12,16;4,3)$:
\[
\begin{footnotesize} \begin{array}{|r|r|r|r|}\hline
I & -I & I & -I  \\ \hline
-I & I & -I & I \\ \hline
I & -I & I & -I \\ \hline
\end{array} 
\end{footnotesize}
\]
\end{ex}

We end this section by constructing NASMs whose row and column weights are even. 

\begin{theorem}\label{thm:nonuniform_2}
Let $\hh = (h_1, h_1, \ldots, h_m, h_m)$ and $\kk=(k_1, k_1, \ldots, k_n, k_n)$ be sequences of positive integers.
  There exists a NASM$(2m, 2n; 2\hh, 2\kk)$ if and only if condition
  \eqref{GHA:nec2} holds.
\end{theorem}
\begin{proof}
  By Remark \ref{NASM:nec}, we only need to show sufficiency. 
  Since $\hh$ and $\kk$ satisfy condition \eqref{GHA:nec2},
  one can check that the same holds for $\hh'=(h_1, \ldots, h_m)$ and 
  $\kk'=(k_1, \ldots, k_n)$. Then, Theorem \ref{GR} guarantees the existence of an $m\times n$ matrix, say $A$, over $\mathbb{Z}_2$ such that $w_r(A)=\hh'$ and $w_c(A)=\kk'$. The array obtained from $A$ by replacing each $0$ with
  $\left(\begin{matrix}
    0 & 0\\
    0 & 0
  \end{matrix}\right)$
  and each $1$ with
    $\left(\begin{matrix}
    1 & -1\\
    -1 & 1
  \end{matrix}\right)$ is the desired NASM.
\end{proof}

Further constructions of NASMs will be provided in a paper in preparation
\cite{Tr}.

\section{Simple GHAs over a cyclic group} 
\label{sec:results_2}
\label{globallysimple}

In this section, we build simple GHAs over a cyclic group.
We start by showing that the alternated form of an increasing sequence of integers is simple modulo $v$ for every sufficiently large $v$.

\begin{lemma}\label{plus_minus_a:lem}
Let $\mathbf{a}=(a_1, \ldots, a_n)$ be an increasing sequence of positive integers, and let $\mathbf{b}\in\{\mathbf{a}^\pm, - \mathbf{a}^\pm\}$ and $v> a_n$. Then, 
all runs of $\mathbf{b}$ are nonzero $(\mathrm{mod}\; v)$.
\end{lemma}
\begin{proof}
By Remark \ref{plus_minus_a:rem}, it is enough to prove the assertion 
when $\mathbf{b} = \mathbf{a}^\pm = (-a_1, a_2, -a_3, \ldots, (-1)^na_n)$. 
We start by showing that the sum $s_0(\bb)$ of $\bb$ is nonzero modulo $v$. 
Set $a_0=0$, and let $x=a_{n}$ if $n$ is odd, otherwise set $x=0$.
We notice that
  $1\leq -a_{2j-1} + a_{2j} \leq -a_{2j-2} + a_{2j} -1,$
  for  $1 \leq j \leq \lfloor n/2 \rfloor$. 
  Therefore,
  \[
  \begin{aligned}
    \lfloor n/2 \rfloor - x \leq s_0(\bb)
     &= \sum_{j=1}^{\lfloor n/2 \rfloor} (-a_{2j-1} + a_{2j}) - x 
        \leq \sum_{j=1}^{\lfloor n/2 \rfloor} (-a_{2j-2} + a_{2j} -1) - x 
        \\
     &= a_{2\lfloor n/2 \rfloor} - x - \lfloor n/2 \rfloor.
  \end{aligned}
  \]
  In other words,
  \[
  \begin{aligned}   
    & \text{if $n$ is even, then $x=0$ and $0 < n/2 \leq 
    s_0(\bb) 
    \leq a_n - n/2$}; \\
    & \text{if $n$ is odd, then $x=a_n$ and $\lfloor n/2\rfloor - a_n \leq 
        s_0(\bb)
    \leq a_{n-1} - a_{n} - \lfloor n/2\rfloor<0$.} 
  \end{aligned}
  \]
Hence, $s_0(\bb) \not\equiv 0 \;(\mathrm{mod}\; v)$. This means that the sums of the alternated forms of an increasing sequence of integers are nonzero modulo $v$, provided that $v$ is larger than the maximum integer in the sequence.

For every $1\leq i<j\leq n$, set $\aa_{ij}=(a_i, a_{i+1}, \ldots, a_{j})$ and 
$\bb_{ij}=\aa_{ij}^{\pm}$. Clearly, $s_{ij}(\bb) = s_0(\bb_{ij})$.
Since $\aa_{ij}$ is increasing and $v>a_n\geq a_{j}$, by the first part of the proof it follows that $s_{ij}(\bb)\not\equiv 0 \;(\mathrm{mod}\; v)$.
\end{proof}

The previous construction can be slightly modified to obtain simple sequences whose total sum is zero.
\begin{lemma}\label{plus_minus_b:lem}
Let $\mathbf{a}=(a_1, \ldots, a_n)$ be a sequence of distinct positive integers, such that $a_1<a_2<\cdots<a_{n-1}$ and 
let $\mathbf{b}\in\{\mathbf{a}^\pm, - \mathbf{a}^\pm\}$ and $v> a_n$.
If $s_0(\mathbf{b}) \equiv 0 \;(\mathrm{mod}\; v)$, then 
all proper runs of $\mathbf{b}$ are nonzero $(\mathrm{mod}\; v)$.
\end{lemma}
\begin{proof} By Remark \ref{plus_minus_a:rem}, it is enough to prove the assertion when $\mathbf{b} = \mathbf{a}^\pm$. 
Since $(a_1, \ldots, a_{n-1})$ is increasing, by Lemma \ref{plus_minus_a:lem} we have that 
all runs of $(-a_1, a_2,$ $-a_3, \ldots, (-1)^{n-1}a_{n-1})$ 
are nonzero modulo $v$. 
Considering that by assumption $s_0(\mathbf{b}) \equiv 0 \;(\mathrm{mod}\; v)$, it follows that all proper runs of $\mathbf{b}$ are nonzero.
\end{proof}

\begin{rem} Sequences satisfying the assumptions of Lemma \ref{plus_minus_b:lem} can be constructed starting from a balanced sequence, defined in \cite[Definition 3.1]{BMT}. More precisely, an increasing sequence of positive integers $\aa'=(a_1, \ldots, a_{2n})$ is called \textit{balanced} if there exists $\tau\in[1,n]$ such that 
$s_0(\bb^{\pm}) = s_0(\cc^{\pm})$, where $\bb=(a_1,\ldots, a_{2\tau})$ and
$\cc=(a_{2\tau+1},\ldots, a_{2n})$. In \cite{BMT}, using a different terminology, it is shown that the sequence
$\aa=(a_1, \ldots, a_{2\tau}, a_{2\tau+2}, \ldots, a_{2n}, a_{2\tau+1})$
and its alternated forms  $\aa^\pm$ and $-\aa^\pm$ satisfy the assumptions of Lemma \ref{plus_minus_b:lem}, hence $\aa^\pm$ and $-\aa^\pm$ 
are zero sum and simple modulo $v> a_{2n}$.
This result is used in the proof of Theorem \ref{fromNASMtoGHA2}.
\end{rem}

Lemmas \ref{plus_minus_a:lem} and \ref{plus_minus_b:lem} are used in the following, together with near alternating sign matrices, to build simple GHAs.

\subsection{Nonzero sum and simple GHAs}
\label{globallysimple:nonzero sum}

\begin{theorem}\label{fromNASMtoGHA} 
Assume there is a NASM$(m,n;\hh,\kk)$. Then, there exists a \nz{} and simple NGHA$_S(m,n;\hh,\kk)$ over $\Z_v$ if and only if the necessary condition \eqref{GHA:nec} holds, that is, $S$ has no involutions and $|S|= s_0(\hh) = s_0(\kk)$.
\end{theorem}
\begin{proof} The necessity of condition \eqref{GHA:nec} has been discussed in the introduction; therefore, it is enough to prove sufficiency.

Let $A$ be a NASM$(m, n; \hh,\kk)$ and $S \subseteq\Z_v^+$ a set satisfying the assumptions. Note that $v \geq 2|S|+1$ and $|S|$ is the weight of $A$.
Furthermore, set $H=\pi_v(A\circ f(A^*))$ where 
$f:[0,|S|]\rightarrow [0, \lfloor \frac{v-1}{2} \rfloor]$ is the increasing map fixing $0$ such that $S=\pi_v(f[1,|S|])$. 

We claim that $H$ is the desired NGHA.
First, we notice that $\|\mathcal{E}(H)\| = \pi_v(f(A^*)) = \pi_v(f[1,|S|]) = S$. 
Now, set $\mathbf{a}_i = f(A^*)_{(i)}$,
$\bb_i = A_{(i)}\circ \mathbf{a}_i$, for some $i=1,\ldots,m$, and note that $H_{(i)} = \pi_v(\bb_i)$.
Since $A$ is a NASM, then $A_{(i)}$ is alternating, 
hence $\bb_i\in \{\mathbf{a}_i^\pm, -\mathbf{a}_i^\pm\}$.
By the definition of $f$, we have that
$\mathbf{a}_i$ is increasing and its maximum entry is less than $v$. Therefore,
by Lemma \ref{plus_minus_a:lem}, it follows that all runs of $\bb_i$ are nonzero modulo $v$, that is, $H_{(i)}$ is a nonzero sum and simple sequence.
Similarly, one can show that each $H^{(j)}$ is nonzero sum and simple. Therefore, $H$ is the desired NGHA.
\end{proof}

By Theorems \ref{thm:ASM} and \ref{thm:nonuniform_2}, conditions 1 and 2 of Theorem \ref{thm:gl_simple_GH} independently guarantee the existence of a NASM$(m, n; \hh, \kk)$. Therefore, Theorem \ref{thm:gl_simple_GH} is a consequence of Theorem \ref{fromNASMtoGHA}.

\subsection{Zero sum and simple GHAs}
\label{globallysimple:zero sum}
\begin{theorem}\label{fromNASMtoGHA2} 
Let  $T= \pi_v(S)\,\cup\,\pi_v(S+x)$ where 
$S\subseteq \mathbb{N}$ and $x,v\in \mathbb{N}$ satisfy the following properties
\begin{enumerate}
  \item $|S| \equiv 0\pmod{4}$,
  \item $S$ is the disjoint union of pairs of consecutive integers,
  \item $\frac{v}{2} - \max S > x > \max S - \min S$.
\end{enumerate}
If there exists a NASM$(m, n; 4\hh, 2\kk)$ of weight $|S|$,
then there exists a zero sum and simple
GHA$_T(2m,n; (4\hh, 4\hh), 4\kk)$ over $\Z_v$.
\end{theorem}
\begin{proof} Let $\hh=(h_1, \ldots, h_m)\in \mathbb{N}^m$ and let 
$A$ be a NASM$(m, n; \hh, 2\kk)$ of weight $|S| = 4(h_1+\cdots+h_m)$.
Since $S$ is the disjoint union of pairs of consecutive integers, 
we can write $S=\bigcup_{i=1}^m S_i$, with $S_i=\{a_{i1}, \ldots, a_{i,4h_i}\}$, such that
\begin{enumerate}[label=$(\roman*)$] 
  \item $a_{ij} < a_{pq}$ whenever $(i,j)<(p,q)$ (according to the lexicographic order on $\mathbb{N}\times \mathbb{N}$), and
  \label{a1}
  \item $a_{i, 2\ell} = a_{i, 2\ell-1} + 1$,\label{a2}
\end{enumerate}
for every $i,j,\ell, p, q$ belonging to the appropriate range of positive integers.

We build the $m\times n$ matrix $B(\epsilon)$, with $\epsilon\in\{0,1\}$, 
as follows: for every $i\in[1,m]$ and
$j\in [1, 4h_i]$,  we replace (following the natural ordering)
the $j$-th nonzero entry of the $i$-th row of $A$ with $b_{\epsilon, i,j}$, where
\begin{equation}\label{beij}
b_{\epsilon, i,j} = \epsilon x+
\begin{cases}
    a_{ij}, & \text{if $j\in [1, 2h_i]$};\\
    a_{i,j+1},   & \text{if $j\in [2h_i+1, 4h_i-1]$},\\
    a_{i,2h_i+1}, & \text{if $j=4h_i$}.
\end{cases}
\end{equation}

We claim that $C=\pi_v
\left(
\begin{array}{r}
    A \circ B(0)\\
   -A \circ B(1)
\end{array}
\right)
$ is the desired GHA.
First notice that $\|\cE(C)\| = \pi_v(\cE(B(0)))\,\cup\,\pi_v(\cE(B(1))) = T$.

We now show that the $i$-th reduced row of $C$, that is,
\[
C_{(i)} = \pm 
\begin{cases}
\pi_v(A_{(i)} \circ  B(0)_{(i)}) & \text{if $1\leq i\leq m$},\\
\pi_v(A_{(i)} \circ  B(1)_{(i)}) & \text{if $m+1\leq i\leq 2m$}.
\end{cases}.
\] 
is a zero sum and simple sequence.
Let $\bb = B(\epsilon)_{(i)} = 
(b_{\epsilon, i,1}, \ldots, b_{\epsilon, i,4h_i})$ 
for some $\epsilon\in\{0,1\}$ and $i\in[1,m]$. 
By \ref{a1}, we have that 
$(b_{\epsilon, i,1}, \ldots, b_{\epsilon, i,4h_i-1})$ is increasing, and by \ref{a2}, one can check that $s_0(\bb^\pm) =0$;
also, $v>b_{\epsilon, i,4h_i}$ (assumption 3). Therefore, by Lemma
\ref{plus_minus_b:lem}, all proper runs of $\bb$ are nonzero modulo $v$.
Since $A_{(i)} \circ  B(\epsilon)_{(i)}\in \{\bb^{\pm}, -\bb^{\pm}\}$, by Remark \ref{plus_minus_a:rem}, we have that the sequence $C_{(i)}$ is zero sum and simple. 

It is left to show that the $j$-th reduced column $C^{(j)}$ of $C$ is 
a zero sum and simple sequence with respect to some ordering.
Letting $B(0)^{(j)}=(b_1, \ldots, b_{2k_j})$, by construction we have that
\begin{align*}
  C^{(j)} = \pm \pi_v&(-b_1, b_2, -b_3, \ldots, b_{2k_j},  \\
         &\;\; b_1+x, -(b_2+x), b_3+x, \ldots, -(b_{2k_j}+x)).
\end{align*}
Setting $\bb=(b_1, \ldots, b_{2k_j}, b_2+x, b_{2k_j}+x, b_1+x)$,
one can check that $s_0(\bb^\pm) = 0$. 
Since $x > \max S - \min S$ (assumption 3) and in view of condition \ref{a1}, 
removing the last term of $\bb$ produces an increasing sequence. Therefore, 
by Lemma \ref{plus_minus_b:lem}, we have that $\pi_v(\bb^\pm)$ and $-\pi_v(\bb^\pm)$ are zero sum and simple sequences. Since one of them is a reordering of $C^{(j)}$, the assertion follows. 
\end{proof}

\begin{rem}\label{row_shift}
  When $4\hh = \underline{n}$ and $2\kk = \underline{m}$, the array $C$ built in Theorem \ref{fromNASMtoGHA2} can be easily rearranged to obtain a zero sum GHA (with the same list of row and column weigths) that is naturally simple. Indeed, it is enough to move to the end the $(m+1)$-th row of $C$.
\end{rem}

In the following, we focus on (the most studied case for classic Heffter arrays, that is) the case where $S= (\Z_{v}\setminus U)^+$ for some subgroup $U$ of $\Z_{v}$.

\begin{cor}\label{fromNASMtoGHA2:cor1} 
  Let $v=(2d+1)u \equiv u \;(\mathrm{mod}\; 16)$, where $u=1$ or 
  $u\equiv 0 \;(\mathrm{mod}\; 4)$, and let $\hh=(h_1, \ldots, h_m)$ be a partition of $\frac{du}{8}$.
  If there is a NASM$(m, n; 4\hh, 2\kk)$,
  then there exists a zero sum and simple 
  GHA$_{(\Z_v\setminus U)^+}(2m,$ $n; (4\hh, 4\hh), 4\kk)$, where
  $U$ is the subgroup of $\Z_v$ of order $u$.
\end{cor}
\begin{proof} Let $T=S\,\cup\, (S+x)$, where $x=\frac{d}{2}$ and
$S=[1,x]$ when $u=1$, otherwise $x = (2d+1)\frac{u}{4}$ and
$S= 
  \left[1, x\right]\setminus
   \left\{i(2d+1)\mid 1\leq i\leq\frac{u}{4}\right\}$.
Note that $S,x$ and $v$ satisfy the assumptions of Theorem \ref{fromNASMtoGHA2}.
Considering that $\pi_v(T)=(\Z_v\setminus U)^+$ and 
$|S|= \frac{du}{2} = 4(h_1+\cdots+h_m)$, 
the result follows by Theorem \ref{fromNASMtoGHA2}.
\end{proof}

The existence of NASMs provided by Theorems \ref{thm:ASM} and \ref{thm:nonuniform_2}, together with 
Corollary \ref{fromNASMtoGHA2:cor1}, implies 
Theorem \ref{fromNASMtoGHA2:cor2}. 
Note that by Remark \ref{row_shift},
when $2m=4k$ and $n=4h$, by shifting down a suitable row of the 
GHA$_{(\Z_v\setminus U)^+}(4k,4h; 4h, 4k)$  built in Theorem \ref{fromNASMtoGHA2:cor2}.(2), we obtain a zero sum GHA -- with the same parameters -- that is simple according to the natural ordering.

\section{\Nz{} GHAs over an arbitrary group} 
\label{algo}
In this section, we describe an algorithm that builds a \nz{} GHA over an arbitrary group $G$, by replacing the 1's of a suitable $m\times n$ matrix $A$ over $\Z_2$ with the elements of a multisubset of $G$. First, we consider the following subsets of the skeleton of $A$:
\[
\begin{aligned}
 \mathcal{R} &= \{(i,j) \mid \ a_{ij} \neq 0 = a_{i,j+1} = \dots =a_{in}\}, \\
 \mathcal{C} &= \{(i,j) \mid \ a_{ij} \neq 0 = a_{i+1,j} = \dots =a_{mj}\}, \\
 \mathcal{L} &= \mathcal{R} \cap \mathcal{C}. 
\end{aligned}
\]
Clearly, $|\mathcal{R}| = m$, $|\mathcal{C}| = n$ and $1\leq |\mathcal{L}| \leq \min(m,n)$. We now consider the simple graph $\Phi  = (V(\Phi ), E(\Phi))$ having $V(\Phi) = \mathcal{R} \cup \mathcal{C}$ as vertex-set, and whose edge-set 
$E(\Phi)$ is defined as follows:
\begin{equation}\label{adjacency}
E(\Phi ) = 
\big\{\{(i,j), (u,v) \}\mid \text{$(i,j) \in \mathcal{R}$, $(u,v) \in \mathcal{C}$ and either $i=u$ or $j=v$.}\big\}
\end{equation}
We start by showing the following.
\begin{lemma} \label{lem:forest}
$\Phi$ is a forest, and each one of its connected components contains exactly one vertex of $\mathcal{L}$.
\end{lemma}
\begin{proof}
We first show that the smallest vertex (with respect to the lexicographic order) 
of a path of $\Phi $, say $P=[x_1, \ldots, x_\ell]$, is necessarily an end-vertex.
Indeed, assume for a contradiction that there is $i\in [2, \ell-1]$  
such that 
\begin{equation}\label{leastvertex}
\text{
$x_i = (u_i, v_i) < x_j = (u_j, v_j)$,\;\;\; for every $j\in [1, \ell]\setminus\{i\}$.
}
\end{equation}
By the definition of the graph $\Phi $, if $x_i\in \cR$, then $x_{i-1}, x_{i+1}\in \cC$. Hence, by \eqref{adjacency} and \eqref{leastvertex}, 
it follows that either $u_i < u_j$ and $v_i = v_j$, or $u_i = u_j$ and $v_i < v_j$, for each $j \in \{i-1,i+1\}$. The latter case cannot happen, since $x_i \in \cR$ is the last cell of its row with a nonzero entry. Therefore, $u_i < u_j$ and $v_i = v_j$
for each $j \in \{i-1,i+1\}$, which however contradicts the fact that $x_{i-1}, x_{i+1}\in \cC$. With a similar reasoning when $x_i\in \cC$, we obtain another contradiction. Therefore, the smallest vertex of a path of $\Phi $ is one of its end-vertices.

This property of $\Phi $ implies that it does not contain cycles, hence $\Phi $ is a forest.
Also, there is no path joining two vertices of $\cL$. Indeed, if $P=[x_1, \ldots, x_\ell]$ is such a path, and $x_1$ is its smallest vertex, then by definition of $\cL$, we would have $x_2<x_1$, contradicting the minimality of $x_1$. Therefore, the vertices of $\cL$ belong to different components of $\Phi $.
\end{proof}

\begin{ex}\label{ex:forest}
Let $m = 6$, $n=8$, 
$\mathbf{h} = (4,3,3,2,2,5)$ and $\mathbf{k}= (3,2,3,2,3,3,1,2)$. We construct the following $6 \times 8$ matrix $A$ with elements in $\Z_2$ such that $w_r(A) = \mathbf{h}$ and $w_c(A)=\mathbf{k}$:
\[
\begin{array}{|r|r|r|r|r|r|r|r|} \hline
1&1 &1 &0 &\bf{1} &0 &0 &0 \\ \hline
1& 0& 1& 0& 0&\bf{1} &0 &0\\ \hline
0&\bf{1} &0 &1 & \bf{1}&0 &0 &0\\ \hline
0&0 &\bf{1} & 0&\bf{1} &0 &0 &0\\ \hline
0& 0& 0&0 & 0& 1&0 &\bf{1} \\ \hline
\bf{1}& 0& 0&\bf{1} &0 &\bf{1} &\bf{1} &\bf{1}\\ \hline
\end{array}
\]
where in bold we have highlighted the cells corresponding to $V(\Phi)$. We then have:
\[
\begin{aligned}
\mathcal{R}&= \{(1,5),(2,6),(3,5),(4,5),(5,8),(6,8)\}, \\
\mathcal{C}&= \{(6,1),(3,2),(4,3),(6,4),(4,5),(6,6),(6,7),(6,8)\}, \\
\mathcal{L}&= \{(4,5),(6,8)\}.
\end{aligned}
\]
We conclude with a drawing of  $\Phi$:
\begin{center}
\includegraphics[width =0.9 \textwidth]{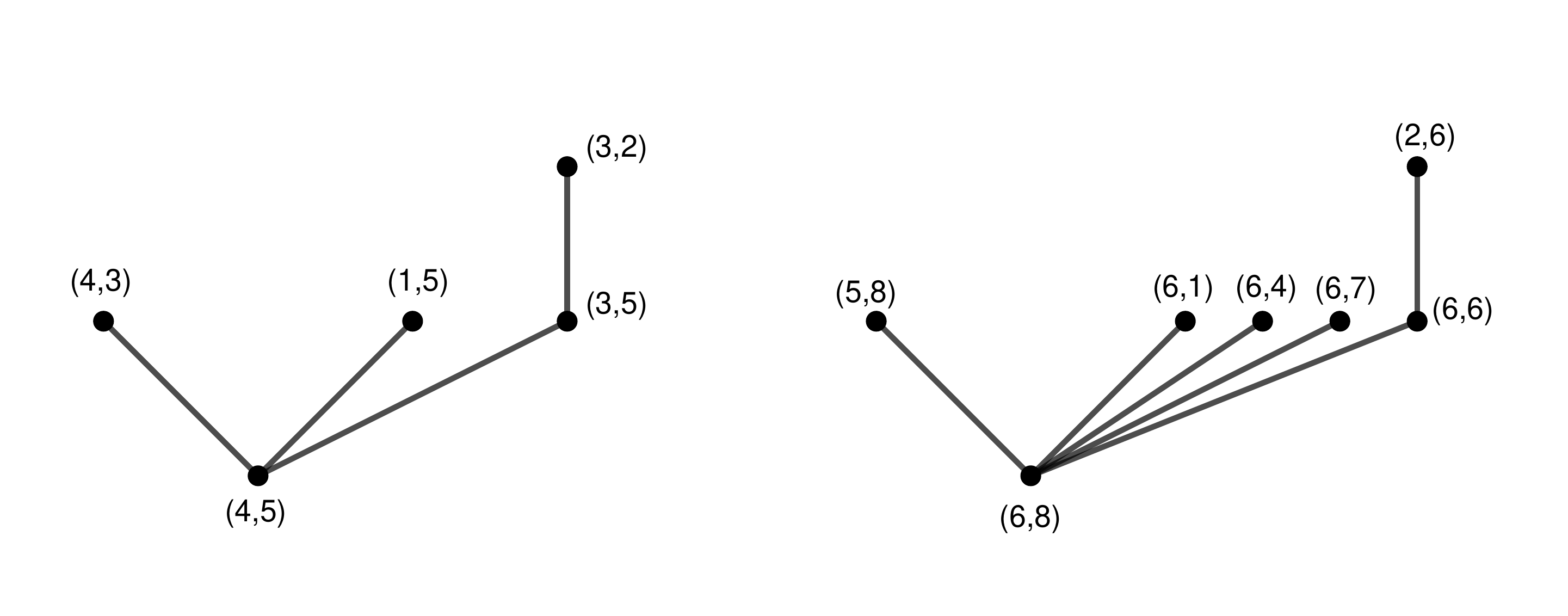}
\end{center}
\end{ex}

We are now ready to prove the main result of this section.
\begin{theorem}\label{thm:main_gh-1}
Let $A$ be an $m \times n$ array over $\Z_2$ such that $w_r(A)=\hh$ and $w_c(A)=\kk$. Let $G$ be an arbitrary group and let $S\subseteq G^+$ such that
\begin{equation}
\label{eq:suff}
 |S\setminus I(S)| \geq |\cR \cup \cC|.
\end{equation}
There exists a \nz{} 
NGHA$^\lambda_S(m, n; \mathbf{h}, \mathbf{k})$ over $G$ if and only if the necessary condition \eqref{GHA:nec} holds.
\end{theorem}

\begin{proof} Let $A = (a_{i,j})$, $\hh$, $\kk$ and $S$ be as in the statement. 
By cyclically shifting the rows (if necessary), we can assume that 
\begin{equation}\label{shift}
  \text{the last row of $A$ has the largest weight}. 
\end{equation}
Let $\mathcal{R}$, $\mathcal{C}$ and $\mathcal{L}$ be the sets of cells of $A$ defined above, and let $\Phi$ be the forest with vertex-set $V(\Phi)=\cR\,\cup\,\cC$, associated to $A$. 
In view of Lemma \ref{lem:forest}, for every $z \in \cL$,
we denote by $\Phi_z$ the connected components of $\Phi$ containing $z$ 
and consider a maximal path $P_z$ in $\Phi_z$ that contains $z$ as an end-vertex. 

Let $T=\,^\lambda(S\setminus I(S))\, \cup\,^\frac{\lambda}{2} I(S)$.
By condition \eqref{eq:suff}, there exists a subset 
$F$ of $S\setminus I(S)$ 
having the same cardinality as $|V(\mathcal{F})|$.
We are going to describe an algorithm that will construct the desired GHA by \emph{properly filling} the cells of $A$, that is, by replacing the 1's in $A$ with all the elements of $T$ in such a way that the rows and columns of the final matrix have nonzero sums. 

\paragraph{1. Filling $\cB = skel(A) \setminus (\cR \, \cup \, \cC)$.} 
We arbitrarily fill the cells of $A$ belonging to 
$\cB$ with $T\setminus F$.

\paragraph{2. Filling isolated vertices of $\Phi$.} 
Let $\cL' \subseteq \cL$ be the set of isolated vertices of $\Phi$.
For every $(i,j)\in \cL'$ and $g\in G$, we say that $g$ is a forbidden element for $(i,j)$ if at least one of the following conditions holds:
\begin{enumerate}
  \item[($a$)] $(i,j)\in \cR$ and by filling the cell $(i,j)$ of $A$ with $g$ we have $\sigma_0(A_i)=0$;
  \item[($b$)] $(i,j)\in \cC$ and by filling the cell $(i,j)$ of $A$ with $g$ we have $\sigma_0(A^j)=0$.
\end{enumerate}
It follows that each cell in $\cL'$ has at most two forbidden elements.
Note that $|\cL'| \leq |\cL| \leq |V(\Phi)|$. If $|\cL'| < |\cL|$, considering that the elements of $F$ are pairwise distinct in absolute value, it is not difficult to check that
there is an injective map $\varphi:\cL' \rightarrow F$ such that either $\varphi(i,j)$ or $-\varphi(i,j)$ is not a forbidden element for $(i,j)$. We use this element to fill the cell $(i,j)$ of $A$, for every $(i,j)\in \cL'$.
 This procedure guarantees that all rows and columns through a cell of $\cL'$ have a nonzero sum.  We then remove $\varphi(\cL')$ from $F$.

Now assume that $|\cL'| = |\cL|$. This means that all components of $\Phi$ consist of isolated vertices, that is, $\cL=\cR=\cC$. Therefore, since all cells in the last row of $A$ belong to $\cC$, it follows that the weight of $A_m$ is 1, and by condition \eqref{shift}, each row of $A$ has weight 1. In this case, $0$ is the only forbidden element for each of the isolated vertices. Hence, we can arbitrarily fill the cells of $A$ belonging to $\cL'$ with the elements of $F$, without creating rows or columns whose sum is $0$.

\paragraph{3. Filling  $W=V(\Phi_z)\setminus V(P_z)$, for every $z \in \cL\setminus \cL'$.} 
As before, for every $(i,j)\in W$, $g\in G$ is a forbidden element for
$(i,j)$ if either condition $(a)$ or $(b)$ (in step 2) is satisfied. 
Note that only one between $(a)$ and $(b)$ holds, since $W\cap \cL=\varnothing$, hence
each cell in $W$ has exactly one forbidden element. Since $|F|>|W|$ and all the elements of $F$ are pairwise distinct, there is
an injective map $\varphi: W \rightarrow F$ such that $\varphi(i,j)$ is none of the forbidden elements for $(i,j)$. This guarantees that 
filling each cell $(i,j)$ of $A$ with 
$\varphi(i,j)$, for every $(i,j)\in W$ does not produce rows or columns whose sum zero.
We then remove $\varphi(W)$ from $F$.

\paragraph{4. Filling $V(P_z)$, for every $z \in \cL\setminus \cL'$.}  
Let $P_z = [x_1, \dotsc, x_\ell=z]$, where each $x_i=(u_i, v_i)$, and let $F_z=\{f_1,\ldots, f_\ell\}$ be any subset of $F$ of size $\ell$.
If $u_{\ell-1} = u_\ell$ (resp. $v_{\ell-1} = v_\ell$), 
let $g_\ell\in G$ such that by filling the cell $x_\ell$ of $A$ with $g_\ell$ we have 
$s_0(A^{v_\ell})=0$ (resp. $s_0(A_{u_\ell})=0$). 
In other words, $g_\ell$ is the element that, if used to fill $x_\ell$, makes equal to $0$ the sum of either the row or the column through $x_\ell$. Since $\ell>1$ and the elements of $F_z=\{f_1,\ldots, f_\ell\}$ are distinct in absolute value, we can apply a permutation to $F_z$ so that 
$g_\ell \notin \{f_\ell, -f_\ell\}$.

Now for every $1\leq i \leq \ell-1$  we proceed as follows.
\begin{enumerate}
	\item Let $g_i\in G$ be the unique element in $G$ such that $(u_i, v_i)$ and $g_i$ satisfy either condition $(a)$ or $(b)$ (in step 2);
	\item Fill cell $x_i$ of $A$ with some element in $\{f_i, -f_i\}\setminus\{g_i\}$;  
\end{enumerate}  
This iteration guarantees that all rows and columns through the cells of $V(P_z)$ have nonzero sums, except possibly for the row and column through $x_{\ell}$. The two forbidden elements for $x_\ell$ are $g_\ell$ and $g'_\ell$. Recall that $g_\ell \notin \{f_\ell, -f_\ell\}$, therefore we can fill $x_\ell$ with some element in 
$\{f_\ell, -f_\ell\}\setminus \{g_\ell , g'_\ell\}$ so that the row and column through 
$x_\ell$ have a nonzero sum. \\

It follows that $A$ is a GHA$^\lambda_S(m, n; \mathbf{h}', \mathbf{k})$ over $G$, where
$\mathbf{h}'$ is a cyclic permutation of $\hh$. Since a cyclic shift of the rows does not change the order of the column sums, by shifting back the rows of $A$ (if necessary), we obtain the desired GHA.
\end{proof}

By using the trivial upper bound  $|\cR \cup \cC|\leq m+n-1$, we obtain Theorem \ref{thm:main_gh}.
The latter improves previous results contained in \cite{CDF}, where the existence of a \nz{} uniform GHA$_S^\lambda $ is proven, using a probabilistic and non-constructive approach, when $|G|\geq 41$ and 
$S = (G\setminus J )^+$, where $J$ is a subgroup of $G$.

\begin{ex}
Let $A$ be the following $6\times 8$ array having $w_r(A) =\mathbf{h}= (5,6,5,4,5,5)$ and $w_c(A) =\mathbf{k}= (6,3,4,6,4,3,2,2)$ (notice that its graph $\Phi$ is isomorphic to the one of the array constructed in Example \ref{ex:forest}):
\[
\begin{array}{|r|r|r|r|r|r|r|r|} \hline
1&1 &1 &1 &\bf{1} &0 &0 &0 \\ \hline
1& 1& 1& 1& 1&\bf{1} &0 &0\\ \hline
1&\bf{1} &1 &1 & \bf{1}&0 &0 &0\\ \hline
1&0 &\bf{1} & 1&\bf{1} &0 &0 &0\\ \hline
1& 0& 0&1 & 0& 1&1 &\bf{1} \\ \hline
\bf{1}& 0& 0&\bf{1} &0 &\bf{1} &\bf{1} &\bf{1}\\ \hline
\end{array}
\]
We have $|\cR \cup \cC| = 12$ and $w(A) = 30$, hence for  $S= \{\alpha^j: j\in[1,12]\} \cup \{ \alpha^j  \beta: j \in[0,5] \}$ subset of the dihedral group 
$Dih_{25} = \langle \alpha, \beta | \alpha^{25} = \beta^2 = (\alpha\beta)^2 =  1 \rangle$, having order $50$, we  construct a 
NGHA$^2_S(6,8;\mathbf{h},\mathbf{k})$ by following the proof of Theorem \ref{thm:main_gh}. We point out to the reader that, in this example, we use the (more standard) multiplicative notation for the dihedral group $Dih_{25}$.

We begin by defining $F = S\setminus I(S) = \{\alpha^j : j\in[1,12]\}$, and we fill the array $B$ with  the elements of 
$(^2S) \setminus (F \cup I(S))=S$:
\[
\begin{array}{|r|r|r|r|r|r|r|r|} \hline
\beta& \alpha \beta & \alpha^2 \beta &\alpha^2&\bf{1} &0 &0 &0 \\ \hline
\alpha^3 \beta& \alpha&\alpha^4 \beta& \alpha^4& \alpha^3&\bf{1} &0 &0\\ \hline
\alpha^7&\bf{1}& \alpha^6 &\alpha^5\beta & \bf{1}&0 &0 &0\\ \hline
\alpha^{10}&0 &\bf{1} & \alpha^5&\bf{1} &0 &0 &0\\ \hline
\alpha^{11}& 0& 0&\alpha^8 & 0&  \alpha^9&\alpha^{12}& \bf{1}\\ \hline
\bf{1}& 0& 0&\bf{1} &0 &\bf{1} &\bf{1} &\bf{1}\\ \hline
\end{array}
\]
The graph $\Phi$ does not contain isolated vertices, thus we directly move to step 3. of the proof of Theorem \ref{thm:main_gh} by filling $V(\Phi) \setminus V(P_1\cup P_2)$, where $P_1 = [(6,8),(6,1)]$ and $P_2 = [(4,5), (3,5),(3,2)]$. 
We can fill each one of these cells with either $a$ or $a^{-1}$, for any $a$ in $F$:
\[
\begin{array}{|r|r|r|r|r|r|r|r|} \hline
\beta& \alpha \beta & \alpha^2 \beta &\alpha^2&\alpha^{-6} &0 &0 &0 \\ \hline
\alpha^3 \beta& \alpha&\alpha^4 \beta& \alpha^4& \alpha^3&\alpha^7 &0 &0\\ \hline
\alpha^7&\bf{1}& \alpha^6 &\alpha^5\beta & \bf{1}&0 &0 &0\\ \hline
\alpha^{10}&0 & \alpha^{4} & \alpha^5&\bf{1} &0 &0 &0\\ \hline\alpha^{11}& 0& 0&\alpha^8 & 0&  \alpha^9&\alpha^{12}&\alpha^2\\ \hline
\bf{1}&  0& 0&\alpha^5 &0 &\alpha &\alpha^3  &\bf{1}\\ \hline
\end{array}
\]
In step 4. we have $F_1 = \{ \alpha^8,\alpha^9\}$ and $F_2 = \{\alpha^{10},\alpha^{11},\alpha^{12}\}$. Since no rows or columns of the array have forbidden elements contained in $F_1 \cup F_2$, we can fill the remaining cells as follows:
\[
\begin{array}{|r|r|r|r|r|r|r|r|} \hline
\beta& \alpha \beta & \alpha^2 \beta &\alpha^2&\alpha^{-6} &0 &0 &0 \\ \hline
\alpha^3 \beta& \alpha&\alpha^4 \beta& \alpha^4& \alpha^3&\alpha^7 &0 &0\\ \hline
\alpha^7&\alpha^{10}& \alpha^6 &\alpha^5\beta &\alpha^{-11}&0 &0 &0\\ \hline
\alpha^{10}&0 & \alpha^{4} & \alpha^5&\alpha^{12} &0 &0 &0\\ \hline\alpha^{11}& 0& 0&\alpha^8 & 0&  \alpha^9&\alpha^{12}&\alpha^2\\ \hline
\alpha^8&  0& 0&\alpha^5 &0 &\alpha &\alpha^3  &\alpha^{9}\\ \hline
\end{array}
\]

\end{ex}

\section{GHAs and orthogonal decompositions}
\label{sec:decs}
Given a graph $\Gamma$, directed or not, possibly with multiple edges, we denote by $V(\Gamma)$ its set of vertices and by $E(\Gamma)$ its multiset of arcs or edges. We speak of a 
\emph{multigraph} if some edges of $\Gamma$ has multiplicity greater than 2, otherwise the graph is \emph{simple}.
We denote by $^\lambda \Gamma$ the \emph{$\lambda$-fold} of a graph $\Gamma$, that is, the graph on the same vertex-set as $\Gamma$ obtained by repeating $\lambda$ times every edge of $\Gamma$. 
If $\Gamma$ is directed, we denote by $\Gamma^*$ the \emph{undirected version} of $\Gamma$, that is, the undirected graph obtained by replacing each arc of $\Gamma$ with an undirected edge.
Note that if $\Gamma$ is directed, then clearly $(^\lambda \Gamma)^* = \,^\lambda (\Gamma^*)$. Therefore, we simply write $^\lambda \Gamma^*$.

Recall that a decomposition of a graph $\Gamma$ is a set $\cF = \{F_1, \ldots, F_n\}$ of subgraphs of
$\Gamma$ whose edge-sets partition between them $E(\Gamma)$. 

\begin{rem}\label{rem:6:1}
Note that if $\cF=\{F_1, \ldots, F_n\}$ is a decomposition of a directed multigraph $\Gamma$, then $\cF^*=\{F^*_1, \ldots, F^*_n\}$ decomposes 
$\Gamma^*$.
\end{rem}

Two decompositions of
$\Gamma$, say $\cF$ and $\cF'$, are said to be \emph{orthogonal} if $|E(F)\,\cap\,E(F')|\leq 1$, for every $F\in\cF$ and $F'\in\cF'$.

Let $G$ be a group and take a set $D\subseteq G\setminus\{0\}$. 
The \emph{Cayley digraph}  $\cayd{G}{D}$ is the simple directed graph with vertex set $G$ whose arcs are exactly those of the form $(d+x, x)$ for every $x\in G$ and $d\in D$. Similarly, the \emph{Cayley graph}
$\cay{G}{D}$ is the simple undirected graph with vertex set $G$ whose edges are exactly those of the form $\{x, d+x\}$ for every $x\in G$ and $d\in D$. In both cases, we refer to $D$ as the connection set. Note that $\cay{G}{D} = \cay{G}{D^+}$. 

Letting $M= \bigcup_{i=1}^n \,^{\lambda_i} D_i$ be a multiset of elements of $G\setminus\{0\}$, where $D_1, \ldots,$ $D_n$ are pairwise disjoint sets, 
we refer to
\[
\cayd{G}{M} = \bigcup_{i=1}^n \,^{\lambda_i} \cayd{G}{D_i}\;\;\;\text{and}\;\;\;
\cay{G}{M} = \bigcup_{i=1}^n \,^{\lambda_i} \cay{G}{D_i},
\]
as the Cayley digraph and the Cayley graph, respectively, with connection multiset $M$. 
Clearly, the above graphs are simple if and only if $M$ is a set.

\begin{rem}\label{rem:6:2}
Note that $\cayd{G}{M}^* = \cayd{G}{\|M\|}^*$.
\end{rem}

Assuming that the directed graph $\Gamma$ is simple, then $\Gamma^*$ has edges of multiplicity either $1$ or $2$; also, an edge $\{x,y\}\in E(\Gamma^*)$, with $x\neq y$, has multiplicity $2$ if and only if $(x,y), (y,x)\in E(\Gamma)$. In particular,
for a given connection set $D$, we have that
\begin{enumerate}
\item $\cayd{G}{D}^*$ is simple if and only if $D\cap -D = \varnothing$, 
      and in this case, $\cayd{G}{D}^* = \cay{G}{D}$;
\item $\cayd{G}{D}^*$ is the $2$-fold of a simple graph if and only if $D=-D$, 
      and in this case, $\cayd{G}{D}^{*} =\,^2\cay{G}{D}$.
\end{enumerate}
One can easily check that $\cayd{G}{M}$ is the $\lambda$-fold of a simple directed graph if and only if $M=\,^\lambda D$ for some set $D\subset G\setminus\{0\}$. Similarly,
$\cayd{G}{M}^*$ is the $\lambda$-fold of a simple undirected graph if and only if 
$\|M\|=\,^{\lambda/2} I(D)\,\cup\,^{\lambda} D'$ for some set $D\subseteq G^+$, where
$I(D)$ is the set of all involutions of $D$ and $D'=D\setminus I(D)$. Indeed, by noticing that $I(D)=-I(D)$ and $D'\,\cap\,-D=\varnothing$, we have that
\begin{align*}
  \cayd{G}{M}^* &= \cayd{G}{\|M\|}^* =
  \,^{\lambda/2} \cayd{G}{I(D)}^* \,\cup\,  \,^{\lambda} \cayd{G}{D'}^* \\
  &= \,^{\lambda} \cay{G}{I(D)} \,\cup\,  \,^{\lambda} \cay{G}{D'}\\
  &= \,^{\lambda} \cay{G}{D}.  
\end{align*}

\subsection{Orthogonal decompositions of Cayley (di)graphs}
Given a graph $\Gamma$, directed or not, with vertices in a group $G$, we denote by 
\begin{enumerate} 
\item $\Gamma+g$ the \emph{right translate} of $\Gamma$ by an element $g\in G$, that is, the graph obtained from $\Gamma$ by replacing each vertex $x$ with $x+g$;
\item $Dev(\Gamma)$ the \emph{development} of $\Gamma$, that is, the multiset of right translates of $\Gamma$ by the elements of $G$. 
\end{enumerate}
\begin{rem}\label{rem:6:3}
If $\Gamma$ is directed, then $(\Gamma+g)^* = \Gamma^*+g$.
\end{rem}
We say that a decomposition $\cF$ of $\cayd{G}{S}$ or $\cay{G}{S}$ is 
\textit{$G$-regular} if
$F+g\in \cF$ for every $F\in \cF$ and $g\in G$.
In this case, the group of all right translations induced by $G$ is an \emph{automorphism group} of $\cF$ with a \emph{sharply transitive action} on the vertex-set.

Given a sequence $\omega=(g_1, \ldots, g_n)$ of nonzero elements of a group $G$, let $\omega^-=(g_n, g_{n-1}, \ldots, g_1)$ be the reverse of $\omega$, and set
$s_{i}^-(\omega)=s_{n-i+1,n}(\omega)$.
We denote by $W_+(\omega)$ and $W_-(\omega)$ the following directed walks of length $k$:
\begin{enumerate}
  \item $W_+(\omega) = (0,-s_1(\omega), \ldots, -s_k(\omega))$,
  \item $W_-(\omega) = (0,s^-_1(\omega), \ldots, s^-_k(\omega))$.  
\end{enumerate}

\begin{rem}\label{rem:6:4}
When $\omega$ is simple, 
then $W_+(\omega)$ and $W_-(\omega)$, as well as their undirected versions, are $k$-cycles (i.e. cycles of length $k$) 
or $k$-paths (i.e. paths of length $k$) according to whether $s_0(\omega)=0$ or 
$s_0(\omega)\neq0$, respectively. 
Furthermore, for every $x\in\{+,-\}$,
\begin{itemize}
  \item $Dev(W_x(\omega))$ is a decomposition of
$\cayd{G}{x\cE(\omega)}$, and
  \item given two sequences $\omega_1$ and $\omega_2$ of elements of $G$,  if $|\cE(\omega_1)\,\cap\, \cE(\omega_1)|\leq 1$, then
$Dev(W_x(\omega_1))$ and $Dev(W_x(\omega_2))$ are orthogonal.
\end{itemize}
\end{rem}

Letting $(\omega_r, \omega_c)$ be an ordering of an $m\times n$ matrix $A$ over a group $G$, for every $x\in\{+,-\}$, set 
\[
\cW_x(\omega_r)=\bigcup_{i=1}^m Dev(W_x(\omega_{r,i}))\;\;\; \text{and}\;\;\; 
\cW_x(\omega_c)=\bigcup_{j=1}^n Dev(W_x(\omega_{c,j})).
\]
In view of the above considerations, we obtain the following result whose proof is left to the reader.

\begin{theorem}\label{A_to_decs}
Let $A$ be an $m\times n$ matrix over a group $G$ and let $(\omega_r, \omega_c)$ be an ordering of $A$.
Then, for every $x\in\{+,-\}$, the following holds. 
\begin{enumerate}
 \item  
   $\cW_x(\omega_r)$ $($resp. $\cW_x(\omega_c))$ is a  $G$-regular decomposition of $\cayd{G}{x\cE(A)}$ into walks whose lengths span 
   $^{|G|}\cE(w_r(A))$ $($resp. $^{|G|}\cE(w_c(A)))$.
 
  \item If $\cE(A)$ is a set, then $\cW_x(\omega_r)$ and  
  $\cW_x(\omega_c)$ are orthogonal.
\end{enumerate}
Furthermore, if $(\omega_r, \omega_c)$ is a simple ordering of $A$, 
then each $W_x(\omega_{r,i})$ $($resp. $W_x(\omega_{c,j}))$ is a directed cycle if $s_0(\omega_{r,i})=0$ $($resp. $s_0(\omega_{c,j})=0)$, otherwise it is a directed path.
\end{theorem}

A result similar to Theorem \ref{A_to_decs} holds for undirected Cayley graphs, as shown in the following theorem, where
\[
\cW_x(\omega_r)^*=\bigcup_{i=1}^m Dev(W_x(\omega_{r,i})^*)\;\;\; \text{and}\;\;\; 
\cW_x(\omega_c)^*=\bigcup_{j=1}^n Dev(W_x(\omega_{c,j})^*).
\]
for $x\in\{+,-\}$. 
Given an array $A=(a_{ij})$, we recall that $\|A\|$ denotes
the matrix  $\|A\|= (\|a_{ij}\|)$.
\begin{theorem}\label{A_to_decs_2}
Let $A$ and $A'$ be $m\times n$ matrices over a group $G$ such that $\|A\| = \|A'\|$.
Also, let $\omega_r$ be an ordering of the rows of $A$, and let 
$\omega'_c$ be an ordering of the columns of $A'$. 
Then, for every $x,y\in\{+,-\}$, the following holds. 
\begin{enumerate}
 \item 
   $\cW_x(\omega_r)^*$ $($resp. $\cW_x(\omega'_c)^*)$ is a $G$-regular decomposition of $\cayd{G}{\cE(\|A\|)}^*$ 
into walks whose lengths span $^{|G|}\cE(w_r(A))$ $($resp. $^{|G|}\cE(w_c(A)))$.

 \item If $\cE(\|A\|)$ is a set, then $\cW_x(\omega_r)^*$ and $\cW_y(\omega'_c)^*$ are orthogonal.
\end{enumerate}
Furthermore, if $\omega_r$ and $\omega'_c$ are simple orderings of the rows of $A$ and of the columns of $A'$, respectively, 
then each $W_x(\omega_{r,i})^*$ $($resp. $W_x(\omega'_{c,j})^*)$ is a  cycle if $s_0(\omega_{r,i})=0$ $($resp. $s_0(\omega'_{c,j})=0)$, otherwise it is a path.
\end{theorem}
\begin{proof} By Remarks \ref{rem:6:1}, \ref{rem:6:2}, \ref{rem:6:3} and Theorem \ref{A_to_decs}, it follows that
$\cW_x(\omega_r)^*$ is a decomposition of $\cayd{G}{x\cE(A)}^*=\cayd{G}{\|x\cE(A)\|}^*=\cayd{G}{\cE(\|A\|)}^*$; also, by construction, 
$\cW_x(\omega_r)^*$ is $G$-regular, for $x\in\{+,-\}$. 
Similarly, one can see that $\cW_x(\omega'_c)^*$ is a $G$-regular decomposition of $\cayd{G}{\cE(\|A'\|)}^* = \cayd{G}{\cE(\|A\|)}^*$,
for $x\in\{+,-\}$. 

Now assume that $\cE(\|A\|)$ is a set and suppose (for a contradiction) that  
$\cW_x(\omega_r)^*$ and $\cW_y(\omega'_c)^*$ are not orthogonal, for some $x,y\in\{+,-\}$.
This means that the walks $W_x(\omega_{r,i})^*+g$ and $W_y(\omega'_{c,j})^*+g'$ have at least two edges in common, for some
$i\in[1,m], j\in[1,n]$ and some $g,g'\in G$. It then follows that $\|\cE(\omega_{r,i})\|\,\cap\,\|\cE(\omega'_{c,j})\|\geq 2$, contradicting
the assumption that $\cE(\|A\|)$ has no multiple elements. 

Finally, if $\omega_r$ and $\omega'_c$ are simple orderings of the rows of $A$ and of the columns of $A'$, respectively, then by Theorem
\ref{A_to_decs}, we have that each $W_x(\omega_{r,i})$ 
$($resp. $W_x(\omega'_{c,j}))$ is a directed cycle if $s_0(\omega_{r,i})=0$ $($resp. $s_0(\omega_{c,j})=0)$, otherwise it is a directed path; clearly, the same holds for
$W_x(\omega_{r,i})^*$ $($resp. $W_x(\omega'_{c,j}))^*$, 
for each $x\in\{+,-\}$.
\end{proof}

\begin{rem}
  If $G$ is abelian, it is not difficult to check that Theorem \ref{A_to_decs_2} continues to hold after replacing 
  each $W_x(\omega_{r,i})$ with $-W_x(\omega_{r,i})$ and 
  each $W_x(\omega_{c,j})$ with $-W_x(\omega_{c,j})$.
\end{rem}

Given a GHA, say $H$, by taking $A=H$ in Theorem \ref{A_to_decs} and $A=A'=H$ in Theorem \ref{A_to_decs_2}, we build orthogonal decompositions of Cayley (di)graphs into walks, paths or cycles, thus generalizing similar results for undirected complete (equipartite) graphs based on classic Heffter arrays (see, for example, \cite{CDFP}). 
More precisely, we obtain the following.

\begin{theorem}\label{GHA_to_decs}
Let $H$ be a GHA${}^{\lambda}_S(m, n; \mathbf{h}, \mathbf{k})$ over $G$, and let $(\omega_r, \omega_c)$ be an ordering of $H$.
Then, for each $x,y\in\{+,-\}$, the following holds. 
\begin{enumerate}
 \item  
   $\cW_x(\omega_r)$ $($resp. $\cW_x(\omega_c))$ is a  $G$-regular decomposition of $\cayd{G}{x\cE(H)}$ into walks whose lengths span 
   $^{|G|}\cE(\hh)$ $($resp. $^{|G|}\cE(\kk))$.
   \item If $\cE(H)$ is a set, then $\cW_x(\omega_r)$ and  $\cW_x(\omega_c)$ are orthogonal.

   \item  $\cW_x(\omega_r)^*$ $($resp. $\cW_x(\omega_c)^*)$ is a 
   $G$-regular decomposition of
   $^\lambda \cay{G}{S}$ into walks whose lengths span $^{|G|}\cE(\mathbf{h})$ $($resp. $^{|G|}\cE(\mathbf{k}))$.  
   \item If $\lambda=1$, 
   then $\cW_x(\omega_r)^*$ and  
        $\cW_y(\omega_c)^*$ are orthogonal.
\end{enumerate}
Furthermore, if $(\omega_r, \omega_c)$ is a simple ordering of $H$, 
then $W_x(\omega_{r,i})$ and $W_x(\omega_{r,i})^*$ 
$($resp. $\cW_x(\omega_{c,j})$ and $\cW_x(\omega_{c,j})^*)$ are cycles if $s_0(\omega_{r,i})=0$ $($resp. $s_0(\omega_{c,j})=0)$, otherwise they are paths.
\end{theorem}

The existence results for  \z{} and simple GHAs contained in Section \ref{globallysimple} allow us to construct, via Theorem \ref{GHA_to_decs}, orthogonal cycle decompositions of Cayley graphs. In particular, by Theorems
\ref{fromNASMtoGHA2:cor2}  and \ref{GHA_to_decs}, we obtain the following.
\begin{theorem}\label{thm:ortodecs_4}
  Let $v=(2d+1)u \equiv u \;(\mathrm{mod}\; 16)$, where $u=1$ or $u\equiv 0 \;(\mathrm{mod}\; 4)$. 
  \begin{enumerate}
    \item If $8hm=4kn=du$, then
    there exists two $(\Z_v$-regular$)$ orthogonal cycle decompositions of $K_{2d+1}[u]$ whose cycle-lengths are $4h$ and $4k$, respectively.
    \item 
    Let $\hh=(h_1, h_1, \ldots, h_m, h_m)$ and
    $\kk=(k_1, k_1, \ldots, k_n, k_n)$ be sequences of positive integers 
    such that $8s_0(\hh)=4s_0(\kk)=du$. If $4\hh$ and $2\hh$ satisfy condition
    \eqref{GHA:nec2}, then there exist two $(\Z_v$-regular$)$ orthogonal cycle decompositions of $K_{2d+1}[u]$ whose cycle-lengths span $^{2v}\cE(4\hh)$ and $^{v}\cE(4\kk)$, respectively.
  \end{enumerate}
\end{theorem}

Similarly, the following result constructs
orthogonal decompositions of Cayley graphs into paths of given lengths.

\begin{theorem}\label{thm:ortodecs_3}
Let $S\subseteq \Z_v^+\setminus\{\frac{v}{2}\}$, and let  
$\hh\in\mathbb{N}^m$ and $\kk\in\mathbb{N}^n$ such that $|S|=s_0(\hh)=s_0(\kk)$.
If $\hh$ and $\kk$ satisfy condition \eqref{GHA:nec2}, 
then there exist two $(\Z_v$-regular$)$ orthogonal path decompositions of $\cay{\Z_v}{S}$ whose path-lengths span 
$^{v}\cE(\hh)$ and $^{v}\cE(\kk)$, respectively.
\end{theorem}
\begin{proof} By Theorem \ref{GR}, there is an $m\times n$ matrix $C$ over $\mathbb{Z}_2$ such that $w_r(C)=\hh$ and
$w_c(C)=\kk$. Let $f:[0,|S|]\rightarrow [0, \lfloor \frac{v-1}{2} \rfloor]$ be the increasing map fixing $0$ such that 
$S=\pi_v(f[1,|S|])$. 
Recalling that $C^*$ denotes the position matrix of $C$ (see Section 2), set $B=\pi_v(f(C^*))$.
We denote by $A$ and $A'$ the arrays built as follows:
\begin{enumerate}
  \item to obtain $A$  we alternate, from left to right, the signs of the nonzero entries in each row of $B$;
  \item to obtain $A'$ we alternate, from top to bottom, the signs of the nonzero entries in each column of $B$.
\end{enumerate}
Note that $\|A\| = \|A'\|$, and $\cE(\|A\|) = \cE(\|A'\|) = \cE(B) = \pi_v(f[1,|S|]) = S$.

Let $(\omega_r, \omega_c)$ be the natural ordering of $A$ starting from the first nonzero element in each row and column.
By Theorem \ref{A_to_decs_2}, it follows that $\cW_+(\omega_r)^*$ and $\cW_+(\omega'_c)^*$ are ($\Z_v$-regular) orthogonal
decompositions of 
$\cayd{G}{\cE(\|A\|)}^* = \cayd{G}{S}^* = \cay{G}{S}$. 
It is left to show that each $W_+(\omega_{r,i})^*$ and each $W_+(\omega'_{c,j})^*$ is a path.
Note that $\omega_{r,i} = A_{(i)}$ is an alternated form of 
$B_{(i)} = \pi_v(f(C^*)_{(i)})$.
Since $f(C^*)_{(i)}$ is increasing, by Lemma \ref{plus_minus_a:lem}, 
it follows that $\omega_{r,i}$ is simple and $s_0(\omega_{r,i})\neq0$. By reasoning in a similar way, 
we have that $\omega'_{c,j}$ is simple and $s_0(\omega'_{c,j})\neq0$. The result then follows again by 
Theorem \ref{A_to_decs_2}.
\end{proof}

We end this section by modifying the two decompositions 
$\cW_x(\omega_r)^*$ and $\cW_x(\omega_c)^*$ of a suitable Cayley graph in order to obtain decompositions of the same Cayley graph into closed walks (i.e. circuits). The latter will play a fundamental role in the next section. From now on, we assume $x\in\{+,-\}$.

Recall that the \emph{period} of an element $g\in G$ is the minimum positive integer $p(g)$ such that $p(g)g=0$. 
From elementary group theory, it is known that there is $T\subseteq G$ such that $G=\{kg + t\mid 1\leq k \leq p(g), t\in T\}$.
Given a sequence $\omega$ of nonzero elements of $\Gamma$ such that $0$ and $g$ are the end vertices of $W(\omega)$, we denote by
\begin{enumerate}
  \item $C_x(\omega)=\bigcup_{k=1}^{p(g)} W_x(\omega)^*+kg$ the closed walk (circuit) obtained by joining all translates 
  of the undirected walk $W_x(\omega)^*$ by multiples of $g$,
  \item $\cC_x(\omega) = \{C_x(\omega) + t\mid t\in T\}$.
\end{enumerate}
Note that $|E(C_x(\omega))| = p(g)|E(W_x(\omega))|$. 
Also, $\cC_x(\omega)$ and $Dev(W_x(\omega)^*)$ partition the same multiset of edges,
that is,  $\cC_x(\omega)$ is a decomposition of $\cayd{G}{\cE(\omega)}^*$ into circuits.

\begin{rem}\label{circuit_decs_1}
When $g=0$, that is, $p(g)=1$, we have that $T=G$.
Hence, $C_x(\omega)=W_x(\omega)^*$ and 
$\cC_x(\omega)=Dev(W_x(\omega)^*)$.  Furthermore, if $\omega$ is simple, then by Remark \ref{rem:6:4} $C_x(\omega)$ is an $\ell$-cycle, where $\ell$ is the length of $\omega$, hence
$\cC_x(\omega)$ is an $\ell$-cycle decomposition of 
$\cayd{G}{\cE(\omega)}^*$.
\end{rem}

As in Theorems \ref{A_to_decs_2} and \ref{GHA_to_decs},
letting $A$ be an $m\times n$ matrix over a group $G$
with ordering $(\omega_r, \omega_c)$, one can check that 
\[\cC_x(\omega_r) = \bigcup_{i=1}^m \cC_x(\omega_{r,i}) \;\;\;\text{and}\;\;\; 
  \cC_y(\omega_c) = \bigcup_{j=1}^n \cC_y(\omega_{c,j})
\]
are $G$-regular decompositions of 
$\Gamma=\cayd{G}{\cE(\|A\|)}^*$ into circuits; furthermore, 
if $A$ is a GHA$^{\lambda}_S(m, n; \mathbf{h}, \mathbf{k})$, then
$\Gamma=\,^\lambda\cay{G}{S}$, for every $x,y\in\{+,-\}$.
\begin{rem}\label{circuit_decs_2}
If $A$ is a \z{} and simple GHA$^{\lambda}_S(m, n; \mathbf{h}, \mathbf{k})$ with respect to the ordering $(\omega_r, \omega_c)$,
then by Theorem \ref{GHA_to_decs} and Remark \ref{circuit_decs_1}, we have that
 $\cC_x(\omega_r)=\cW_x(\omega)^*$ and 
$\cC_x(\omega_c)=\cW_x(\omega_c)^*$ 
are cycle decompositions of 
$\,^\lambda \cay{G}{S}$.
\end{rem}

\section{GHAs and biembeddings of graphs onto surfaces}
\label{sec:biembeds}
In this section, we briefly recall the basic definitions 
(see \cite{GrTu, MoTh})
concerning biembeddings of graphs onto surfaces, and we  generalize the constructions of these embeddings based on Heffter arrays to the more general setting of Definition \ref{GHA}.

\begin{defn}\label{defi:embed}
An \textit{embedding} of a graph $\Gamma$ in a surface $\Sigma$ is a continuous injective map $\phi: \Gamma \rightarrow \Sigma$, where $\Gamma$ is viewed as a simplicial $1$-complex
endowed with the usual topology.
\end{defn}
We remark that in our setting $\Sigma$ and $\Gamma$ may not be connected.

A connected component of $\Sigma \setminus \phi(\Gamma)$ is called a $\phi$-\emph{face}. An embedding $\phi$ is said to be \textit{cellular} if every $\phi$-face is homeomorphic to an open disc.
\begin{defn} Let $\cC$ and $\cC'$  be two circuit decompositions of a simple graph $\Gamma$. An embedding $\phi: \Gamma \rightarrow \Sigma$
is called a \emph{biembedding} of $\cC$ and $\cC'$ if 
it is face $2$-colorable and the sets of boundaries of each color class are 
$\phi(\cC)$ and $\phi(\cC')$, respectively. 

\end{defn}

For every edge $e$ of a given graph $\Gamma$, we consider its two possible directions $e^+$ and $e^-$, and we denote by $\tau$ the involution that swaps $e^+$ and $e^-$. For every vertex $v$ in $\Gamma$, a \textit{local rotation} $\rho_v$ is a cyclic permutation of the edges directed out of $v$. If we choose a local rotation for each vertex of $\Gamma$, we then obtain a rotation of the directed edges of $\Gamma$. We also recall the following result (see \cite{GrTu}).

\begin{theorem}\label{thm:graph_rot_eq_embed}
A rotation $\rho$ on $\Gamma$ is equivalent to a cellular embedding of $\Gamma$ in an orientable surface. The face boundaries of the embedding corresponding to $\rho$ are the orbits of $\rho\tau$. 
\end{theorem}

Letting $\Sigma$ be the surface of the embedding built by Theorem
\ref{thm:graph_rot_eq_embed}, we notice that the number of its connected components is equal to the number of connected components of $\Gamma$.
By knowing the number of faces, the genus $g$ of each connected component $C$ of the surface can be obtained from the Euler's formula $v - e + f = 2-2g$,
where $v$, $e$, and $f$ denote the number of vertices, edges and faces determined by the embedding in $C$, respectively.

Now, let $A$ be an $m \times n$ array whose entries are pairwise distinct, that is, $\cE(A)$ does not contain multiple elements and let 
$(\omega_r, \omega_c)$ be an ordering of $A$. 
By abuse of notation, we consider
each $\omega_{r,i}$ (resp. $\omega_{c,j}$) as a permutation of 
$\cE(A_i)$ (resp. $\cE(A^j)$)
and let $\omega_r=\omega_{r,1}\cdots\omega_{r,m}$ and
$\omega_c=\omega_{c,1}\cdots\omega_{c,m}$ denote the two permutations of $\cE(A)$ induced by the orderings $\omega_r, \omega_c$ of the rows and columns of $A$, respectively.
 The orderings  $\omega_r$ and  $\omega_c$ are said to be
\emph{compatible} if $\omega_c  \omega_r$ is a cyclic permutation of length  $|\cE(A)|$. In this case, we refer to $(\omega_r, \omega_c )$ as a compatible ordering of $A$.

We are now able to show how GHAs produce biembeddings of 
Cayley graphs onto orientable surfaces.

\begin{theorem} \label{prop:biembed_cay}
	Let $H$ be a GHA$_S(m,n;\mathbf{h}, \mathbf{k})$ over a group $G$, with a compatible ordering $(\omega_r, \omega_c)$.
	Then there exists a biembedding of the circuit decompositions $\cC_+(\omega_r)$ and 
$\cC_-(\omega_c^-)$ of  $\cay{G}{S}$ over an orientable surface.
\end{theorem}
\begin{proof}
We define a permutation $\rho_0$ of $\pm S = S\,\cup\, -S$ as follows:
\[
\rho_0 = 
\begin{cases}
-\omega_r(a) & \text{if $a \in \cE(H)$,}\\
\omega_c(-a) & \text{if $a \in -\cE(H)$.}
\end{cases}
\]
Notice that if $a\in\cE(H)$, then $\rho_0^2(a) = \omega_c \omega_r (a)$. Since the two orderings $\omega_r$ and $\omega_c$ are compatible, then $\rho_0^2$ acts ciclically on $\cE(H)$. Furthermore, 
for $a \in \cE(H)$ (resp. $a \in -\cE(H)$),  
it can be seen that $\rho_0(a) \in -\cE(H)$ (resp. $\rho_0(a) \in \cE(H)$); hence, $\rho_0$ acts cyclically on $\pm \cE(H) = \pm S$.

Let $\Gamma=\cay{G}{S}$ and define the map $\rho$ on the set $D(\Gamma)$ of directed edges of $\Gamma$ as follows:
\[
\rho((g, a+g)) = (g, \rho_0(a) + g), 
\]
for every $g\in G$ and $a\in S$.
Now, since $\rho_0(a)$ is a cyclic permutation of $\pm S$, 
for each $g \in G$ the map $\rho$ is a local rotation of the edges directed out of $g$; hence, as $g$ ranges through $G$, we have that $\rho$ is a rotation of $\Gamma$. Therefore, by Theorem \ref{thm:graph_rot_eq_embed}, there exists a cellular embedding $\Pi$ of $\Gamma$ in an orientable surface such that the face boundaries are the orbits of $\rho\tau$, where $\tau$ is the involution swapping the direction of every edge, i.e., $\tau((g_1,g_2)) = (g_2,g_1)$ for every $g_1,g_2\in G$.

We now verify that the obtained embedding is comprised of the circuits belonging to $\cC_+(\omega_r)$ and 
$\cC_-(\omega_c^-)$.
Let $a \in \cE(H^j)$ be a nonzero entry of the $j$-th column of $H$ (whose weight is $k_j$) and denote by $p_j$ the order of the sum 
$H^j$ with respect to the ordering $\omega^-_c$. 
We have that:
\begin{align*}
\rho \tau ((g,a + g)) &= \rho(a+g, -a + (a+g)) \\
&=(a+g, \rho_0(-a) +a + g ) \\
&=(a+g, \omega_c(a)+ a+g).
\end{align*}
Hence, $(g, a+g)$
belongs to a face boundary represented by the following circuit of length $\ell=p_jk_j$:
\[
\begin{aligned}
&B=(g, a+g, \omega_c(a)+a+g, \omega^2_c(a) + \omega_c(a)+a+g, \dotsc \\ 
& \dotsc, \omega_c^{\ell-2}(a) + \cdots + \omega^2_c(a) + \omega_c(a)+a+g).
\end{aligned}
\]
Assume that $\omega_{c,j}=(a_1, \ldots, a_{k_j})$ with $a_u=a$
for some $u\in[1, k_j]$. Letting $b = a_{u-1} + a_{u-2} + \cdots + a_1$, one can check that
$B = C_-(\omega_{c,j}^-) -b+g$, therefore $B\in 
\cC_-(\omega_c^-)$.

Consider now the oriented edge $(g,-a+g)$, where $a$ is a nonzero entry of the row $H_i$.
Denote by $q_i$ the order of the sum of $H_i$, with respect to the ordering $\omega_r$.
Since $a \in \cE(A)$, we have that:
\begin{align*}
\rho \tau ((g, -a + g)) &= \rho(-a + g, a + (-a + g)) \\
&= \rho(-a + g, \rho_0(a) + (-a + g))\\
&= (-a + g, -\omega_r(a) -a + g).
\end{align*}
Hence, $(g, -a+g)$ belongs to a face boundary represented by the following circuit of length $\ell=q_ih_i$:
\[
\begin{aligned}
B=
(g, -a+g, -\omega_r(a) -a + g, \ldots, -\omega^{\ell-2}_r(a) \cdots -\omega_r(a) -a + g).
\end{aligned}
\]
Assume that $\omega_{r,i} = (a_1, \ldots, a_{h_i})$ with $a_u=a$ for some $u\in[1,h_i]$. 
Letting $b = -a_{u-1} - a_{u-2} - \cdots - a_1$, one can check that 
$B = C_+(\omega_{r,i}) +b + g$, hence 
$B\in \cC_+(\omega_r)$.

We have therefore shown that every non-oriented edge of $\Gamma$ belongs to the boundaries of exactly two faces of the cellular embedding $\Pi$; these two face boundaries are represented by circuits belonging to $\cC_+(\omega_r)$ and $\cC_-(\omega_c^-)$. Therefore,
the embedding is $2$-colorable and the face boundaries are the images through $\Pi$ of the circuits in $\cC_+(\omega_r)$ and 
$\cC_-(\omega_c^-)$.
\end{proof}

\begin{rem} 
If the GHA in Theorem \ref{prop:biembed_cay} is, in addition, 
\z{} and simple according to the ordering $(\omega_r, \omega_c^-)$,
then by Remark \ref{circuit_decs_2} (the two color classes of face boundaries represented by)
$\cC_+(\omega_r)$ and 
$\cC_-(\omega_c^-)$ are cycle decompositions.
Note that when $G$ is abelian, the property of being \z{} is independent of the chosen ordering, and $(\omega_r, \omega_c^-)$ is simple if and only if $(\omega_r, \omega_c)$ is simple.
\end{rem}

Theorem \ref{prop:biembed_cay} generalizes similar results, in the non abelian case, previously obtained in \cite{A, CDFP, MP23} and concerning classic (relative) Heffter arrays.
Similarly to what has been remarked in \cite{CDFP}, the proof of Theorem \ref{prop:biembed_cay} can be adapted to GHAs with multiplicity greater than 1, thus obtaining the following: 
if $H$ is a GHA$_S^\lambda(m,n;\mathbf{h}, \mathbf{k})$ with a 
compatible ordering $(\omega_r, \omega_c)$, then there exists a biembedding of $^\lambda \cay{G}{S}$ over an orientable surface.

\section{Conclusions and further remarks}
We have introduced and studied the existence of generalized Heffter arrays, briefly GHA$_S^\lambda$ (Definition \ref{GHA}), which extend classic (\nz) Heffter arrays (see the survey \cite{PD23}) by allowing that:
\begin{itemize}
	\item the number of nonzero entries in each row (resp. column) of the array is not constant;
	\item the entries of the GHA, in absolute value, belong to an arbitrary subset $S$ of a group $G$, not necessarily abelian.
\end{itemize}
GHAs can be used to construct orthogonal path or cycle decompositions and biembeddings of Cayley graphs onto orientable surfaces (Sections \ref{sec:decs} and \ref{sec:biembeds}) whose structural properties depend on the sum of the entries in each row and column of a GHA, with respect to a given ordering. We have therefore denoted by
NGHA a generalized Heffter array whose rows and column are ordered according to the natural order: from left to right for the rows, from top to bottom for the columns.

In Section \ref{NASM}, we have introduced and built near alternating sign matrices which we then use in Section \ref{globallysimple} to construct \nz{} and simple NGHAs over $\Z_v$ (Theorem \ref{thm:gl_simple_GH}),
thus widely generalizing previous existence results contained in \cite{CDFP, MP23}. 
In particular, Theorem \ref{thm:gl_simple_GH} settles, in the uniform case, the existence problem for
a \nz{} and simple NGHA$_S$ over $\Z_v$, while Theorem \ref{fromNASMtoGHA2:cor2} 
constructs zero sum and simple GHAs whose row and column weights are congruent to $0$ modulo 4. In the tight case, that is, when row and column weights coincide with the number of columns and rows, respectively, then the classic Heffter array built in Theorem \ref{fromNASMtoGHA2:cor2}.(2) can be rearranged to be naturally simple.
In Section \ref{algo}, we provided an algorithm (Theorem \ref{thm:main_gh-1}) that constructs \nz{} 
NGHAs over an arbitrary group whenever $S$ contains enough noninvolutions. Theorem \ref{thm:main_gh}
extends an existence result 
contained in \cite{CDF}, which is however nonconstructive.

GHAs (and more generally matrices) over an arbitrary group $G$ are used in Section  \ref{sec:decs} to build 
$G$-regular orthogonal decompositions of Cayley graphs, directed or not, into walks, cycles or paths, of various lengths which depend on the row and column weights of $A$. Finally, in Section \ref{sec:biembeds},
we construct biembeddings that originate from suitable GHAs, showing that the same conclusions that hold when starting from classic Heffter arrays remain true for GHAs. 

Uniform multigraphs, that is, graphs whose edges have the same multiplicity $\lambda$, are usually more studied 
since, in many cases, they formalize uniform interactions between objects.
Condition 1 defining a GHA$_S^\lambda$ over a group $G$, say $H$, that is,
\[\;^{2}\|\cE(H)\| = 
   \;^{2\lambda} (S\setminus I(S)) \,\cup\, ^{\lambda}I(S)
\]
is motivated by applications in constructing decompositions and biembeddings
of the Cayley graph $^\lambda \cay{G}{S}$ in which all the edges have the same multiplicity $\lambda$.
However, one may consider further generalizations of these combinatorial matrices with applications to nonuniform graphs.
Indeed, the hypothesis that $\|\cE(H)\|$ covers exactly $\lambda$ times every noninvolution of $S$, and exactly $\lambda/2$ times every involution of $S$, was not used in many cases: see, for example, Theorems \ref{thm:main_gh},
\ref{A_to_decs}, \ref{A_to_decs_2} and \ref{prop:biembed_cay}, thus suggesting that they could remain true in a more general setting.

We could then consider a GHA$^{\boldsymbol{\lambda}}_{\mathbf{S}}(m, n; \mathbf{h}, \mathbf{k})$, with 
$\boldsymbol{\lambda} = (\lambda_1,\dotsc,\lambda_\ell)$
and 
$\boldsymbol{S} = (s_1,\dotsc,s_\ell)\in (G^+)^\ell$,
where now condition 1. of Definition \ref{GHA} becomes:
\[
\| \cE (A)\| = \{^{\lambda_i} s_i \mid 1\leq i\leq\ell \}.
\]
Then, such an array would encode the construction of a walk decomposition of the union of Cayley multigraphs, and if two compatible orderings exist, of its relative biembedding. 
Another interesting variant to study can be obtained by replacing the first condition in Definition \ref{GHA} with 
$\cE (A) = \,^\lambda S$ where $S\subset G\setminus\{0\}$. Indeed, this type of matrices yield, when $\lambda=1$,
pairs of orthogonal decompositions of the directed Cayley graph $\cayd{G}{S}$ (see, Theorem \ref{A_to_decs}).

We ended Section \ref{sec:biembeds} by providing two circuit decompositions, 
denoted by $\cC_x(\omega_r)$ and $\cC_x(\omega_c)$ ($x\in\{+,-\}$),
of a suitable Cayley graph which originate from the rows and columns of an ordered $m\times n$ GHA, say $(H,\omega)$. We pointed out that the lengths of their circuits depend not only on the weights of the rows and columns of $H$, but they also depend on the periods of the sums of each row and column, according to $\omega$. 
We refer to the sequence 
\[
per(H)=(p(\omega_{r,1}), \ldots, p(\omega_{r,m}), p(\omega_{c,1}), \ldots, p(\omega_{c,n}))
\]
of such periods as the \emph{period} of $H$. Note that $per(H)=\underline{1}$ if and only if $H$ is a \z{} GHA.
Also, if $H$ is a GHA over $\Z_p$, with $p$ prime, then $H$ is \nz{} if and only if $per(H)=\underline{p}$.
In future work, we aim to construct uniform GHAs with a constant period that satisfy the following strong property: 
the circuits in $\cC_x(\omega_r)$ and $\cC_x(\omega_c)$ are cycles, that is, 
both 
$\cC_x(\omega_r)$ and $\cC_x(\omega_c)$ are cycle decompositions
of a suitable Cayley graph.

\section*{Acknowledgements}
The authors were partially supported by INdAM-GNSAGA.

\end{document}